\documentclass{article}

\usepackage{cmd}
\usepackage{format}

\title{Building Models of Determinacy from Below}
\author{Obrad Kasum\footnote{Mr.\ Kasum has received funding from the European Union’s Horizon 2020 research and innovation program under the Marie Skłodowska-Curie grant agreement No.\ 945322. \includegraphics[height=2.5mm]{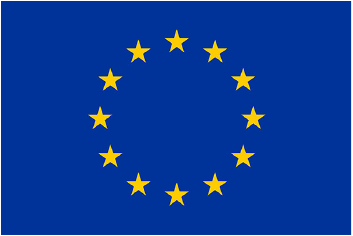}}, Grigor Sargsyan\footnote{Mr.\ Sargsyan's work is funded by the National Science Centre, Poland under the Weave-UNISONO call in the Weave programme, registration number UMO-2021/03/Y/ST1/00281.}}

\begin{document}

\maketitle
\begin{abstract}
    We present an $L$-like construction that produces the minimal model of $\AD_\R+$``$\Theta$ is regular''.
    In fact, our construction can produce any model of $\AD^++\AD_\R+V=L(\ps(\R))$ in which there is no hod mouse with a measurable limit of Woodin cardinals.
\end{abstract}
\tableofcontents
\newpage

\section{Introduction}



\subsection{Statement of the Result}
\label{Statement of the Result}

It was known through the work of Sargsyan that if there exists an inner model containing all reals and satisfying $\AD_\R+$``$\Theta$ is regular'', then there exists the smallest such model.\footnote{
It was shown in \cite{sargsyan2015hod} that if there are two divergent inner models containing all reals and satisfying $\AD^+$, then there exists an inner model containing all reals and satisfying $\AD_\R+$``$\Theta$ is regular'' contained in both.
}
However, this was an existence result rather than a construction, so it was left open what this model looks like.
In this paper we will describe an $L$-like construction that can produce this model, and somewhat larger models.

The minimal model of $\AD_\R+$``$\Theta$ is regular'' satisfies $V=L(\ps(\R))$, so we will focus only on such models.
Each such model $\Mod$ is completely determined by the collection of its sets of reals.
More precisely, if $\Delta$ denotes $\ps(\R)\cap \Mod$, then $\Mod=L(\Delta)$.
So the problem is to construct the powerset of reals of the minimal model of $\AD_\R+$``$\Theta$ is regular''.

Let us call $\Delta\subseteq \ps(\R)$ a \textit{determinacy pointclass} if
$$\ps(\R) \cap L(\R, \Delta) = \Delta$$
and $L(\Delta)\models\AD^+$.
For any set $X$, we denote by $\mu_X$ the club filter on $[X]^\omega$.
Now, given a determinacy pointclass $\Delta$ and an ordinal $\kappa$, we let $\mathscr{O}(\Delta,\kappa)$ be the pointclass
$$\ps(\R) \cap L(\Delta^\omega, \kappa^\omega)[\mu_\Delta, \mu_{\kappa}].$$
For any determinacy pointlcass $\Delta$, if there exists $\kappa$ such that
$$\Delta \subsetneq \mathscr{O}(\Delta, \kappa),$$
we let $\kappa_\Delta$ be the smallest such $\kappa$, otherwise $\kappa_\Delta$ is undefined.
Finally, the pointclass $\mathscr{O}(\Delta)$ is defined if and only if $\kappa_\Delta$ is defined, and it is set to be $\mathscr{O}(\Delta, \kappa_\Delta)$.

The sequence
$$(\De_\alpha : \alpha < \sh)$$
of determinacy pointclasses is defined recursively by iterating the operator $\mathscr{O}$ as long as possible.
If $L(\R)$ does not satisfy $\AD^+$, the sequence is empty, otherwise $\De_0$ is set to be the powerset of reals of $L(\R)$.
If $\De_\alpha$ is defined, then $\De_{\alpha+1}$ is set to be $\mathscr{O}(\De_\alpha)$, if the latter is defined and a determinacy pointclass.
The recursion terminates if the last condition is not met.
At the limit stages we take the constructive closure of the union of the previous pointclasses.
Again, the construction terminates if the obtained pointclass is not a determinacy pointclass.
Our main theorem is then the following.

\begin{mainTheorem}
    Suppose $\Mod$ is an inner model containing all reals and satisfying $\AD^+ + \AD_\R+$``there is no least branch hod pair with a measurable limit of Woodin cardinals''.
    Then there exists $\alpha<\sh$ such that $\ps(\R) \cap \Mod = \De_\alpha.$
    \qed
\end{mainTheorem}

\begin{corollary*}
    If there exists an inner model containing all reals and satisfying $\AD_\R+$``$\Theta$ is regular'', then there exists $\alpha<\sh$ such that $L(\De_\alpha)$ is the minimal model of $\AD_\R+$``$\Theta$ is regular''.
\end{corollary*}
\begin{proof}
    The point is that in the minimal model of $\AD_\R+$``$\Theta$ is regular'', there is no sharp for the minimal hod mouse with an inaccessible limit of Woodin cardinals.
    To see this, let $V$ be this minimal model and suppose $M$ is such a sharp.
    Let $\lambda$ denote the inaccessible limit of Woodin cardinals in it.
    First, the top measure of $M$ can be iterated $\Ord$ times to produce an inner model $M^*$.
    Second, $M^*$ is $\omega_1$-iterable below $\lambda$ by essentially the same strategy as $M$.
    We can then go to the forcing extension by $\Col(\omega,\R)$ and make all reals of $V$ generic over Woodin cardinals of a countable iterate $M^{**}$ of $M^*$.
    The derived model of $M^{**}$ at $\lambda$ is then of the form $L(\R^V, \Delta)$, where $\Delta$ is its powerset of reals, and it satisfies $\AD_\R+$``$\Theta$ is regular'' by \cite[Theorem 1.3]{gappo2022derived}.
    It is easy to see that in fact $\Delta$ is in $V$ and is a surjective image of the reals there, which contradicts the minimality assumption on $V$.
\end{proof}

We will now try to motivate why the operator $\mathscr{O}$ looks the way it does.
The idea behind it is that we want a way to go from any determinacy pointclass to a strictly bigger one.

If there exists an inner model containing all reals and satisfying $\AD^+$, then $L(\R)$ is the least such model.
This minimal model is obtained by a well-known recursive construction:
\begin{eqnarray*}
    J_1(\R) & = & V_{\omega+1} \\
    J_{\alpha+1}(\R) & = & \mathsf{rud}(J_\alpha(\R) \cup \{J_\alpha(\R)\}) \\
    J_\gamma(\R) & = & \bigcup_{\alpha<\gamma} J_\alpha(\R)\quad (\gamma\mbox{ limit}) \\
    L(\R) & = & \bigcup_{\alpha\in\Ord} J_\alpha(\R).
\end{eqnarray*}
In other words, we start with the reals and iterate the operator $J_1$.
From time to time a new set of reals is added to $J_\alpha(\R)$ and by the time we reach $\alpha=\Theta^{L(\R)}$, all sets of reals of $L(\R)$ have been added.
This results in the pointclass $\De_0$ of all sets of reals of $L(\R)$, the $\subseteq$-minimal determinacy pointclass.

A strictly bigger determinacy pointclass cannot be obtain by the simple $L$-construction: an oracle must be added to this purely recursive approach.\footnote{
If there exists a model of $\AD^+$ which extends $L(\R)$ and has a larger $\Theta$, then by \cite[Theorem 10.1.8]{larson2023extensions} this model contains $\R^\sharp$.
In particular, whatever operation is used to get to this model must be able to produce the sharp of the reals.
}
It turns out that if there is a determinacy pointclass strictly larger than $\De_0$, then
$$\ps(\R) \cap L(\De_0^\omega)$$
is such a pointclass.
More generally, the same is true if $\De_0$ is replaced by any determinacy pointclass $\Delta$ such that $\Theta^{L(\Delta)}$ has countable cofinality.\footnote{
The equality $\cof(\Theta^{L(\R)})=\omega$ holds since it is witnessed by the direct limit embedding coming from the HOD analysis for $L(\R)$.
}
So the oracle we use in this case is the sequences of sets in the pointclass.
On the other hand, if $\Theta^{L(\Delta)}$ has uncountable cofinality, then the pointclass
$$\ps(\R) \cap L(\Delta^\omega)$$
is just $\Delta$.\footnote{This follows from the fact that if $\Mod$ is an inner model containing all reals and satisfying $\AD^+$, then for all $\alpha < \Theta^\Mod$, $\alpha^\omega \subseteq \Mod$.}
To continue, we need to distinguish two cases, depending on whether $L(\Delta) \models \AD_\R$.

An \textit{$\AD_\R$ pointclass} is a determinacy pointclass whose constructible closure satisfies $\AD_\R$.
Suppose there exists an $\AD_\R$ pointclass strictly larger than $\Delta$.
In the first case, if $\Delta$ is not an $\AD_\R$ pointclass, we have that
$$\ps(\R) \cap L(\Delta)[\mu_\Delta]$$
is a determinacy pointclass strictly larger than $\Delta$.
On the other hand, in the second case, where $\Delta$ is an $\AD_\R$ pointclass, the answer to how to get a larger determinacy pointclass is more subtle and only partially known.

We would like to say that there exists $\kappa$ such that the powerset of reals of $L(\Delta, \kappa^\omega)[\mu_\kappa]$ is strictly larger than $\Delta$.
However, this can fail in general.\footnote{
For example, it fails if there is a proper class of Woodin limit of Woodin cardinals and if $\Delta$ is the powerset of reals of the Chang-plus model.
We give more details below.
}
What we have been able to show is that there is such a $\kappa$ if $\Delta$ is not too complicated.\footnote{
The upper bound on the complexity of the pointclass is given in terms of strategies for hod mice.
The exact statement is given below.
}
Furthermore, for the smallest $\kappa$ as above, we know that $L(\Delta, \kappa^\omega)[\mu_\kappa]$ satisfies $\AD^+$.
So in this last case, 
$$\ps(\R) \cap L(\Delta, \kappa^\omega)[\mu_\kappa]$$
is a determinacy pointclass strictly larger than $\Delta$.

\subsection{Outline of the Proof}

Before starting the proof of the Main Theorem, we include a section summarizing some facts regarding HOD mice and HOD analysis (Section \ref{HOD Mice and HOD Analysis}).
Near the end of Section \ref{Club Filter}, we start working towards the proof of the Main Theorem.
The introductory part of that section recalls and systemizes some facts concerning the club filter on the collection of countable subsets of a set.
More precisely, we describe the notion of quasi-club filter in the context independent of the Axiom of Choice, which reduces the usual club filter if the Choice holds.

Let now $\Mod$ be an inner model containing all reals and satisfying $\AD^+ + \AD_\R+$``there is no least branch hod pair with a measurable limit of Woodin cardinals''.
We want to show that there exists $\alpha < \sh$ such that
$$\ps(\R) \cap \Mod = \De_\alpha.$$
The proof comes down to showing that the construction does not break down too soon and that a long enough initial segment of the construction is absolute between $V$ and $\Mod$.
The point here is that the sequence of pointclasses that the construction produces is strictly increasing, so as long as it stays in $\Mod$ and does not break down, it will eventually reach the powerset of reals in $\Mod$.

The absoluteness of the construction comes down to the absoluteness of the operator
$$\mathscr{O} : \Delta \mapsto \ps(\R) \cap L(\Delta^\omega, \kappa_\Delta^\omega)[\mu_\Delta, \mu_{\kappa_\Delta}].$$
As long as we know that $\Theta^{L(\Delta)}$ and $\kappa_\Delta$ are strictly smaller than $\Theta^\Mod$, this absoluteness follows from the fact that $\Mod$ contains all the $\omega$-sequences which are bounded in $\Theta$ and from the fact that quasi-club filters are sufficienlty absolute between $\Mod$ and $V$ (the latter is established in Section \ref{Club Filter}).
Thus, the hard part is to show that $\kappa_\Delta$ is well-defined and strictly smaller than $\Theta^\Mod$.

If $\cof(\Theta^{L(\Delta)}) = \omega$, then the joint of any $\omega$-sequence which is Wadge-cofinal in $\Delta$ is a set of reals which is not in $\Delta$.
This shows that $L(\Delta^\omega)$ contains a new set of reals, so $\kappa_\Delta = 0$.
If $\cof(\Theta^{L(\Delta)}) > \omega$, we have to consider two subcases.
First, if $\Delta$ is \underline{not} an $\AD_\R$ pointclass, then $L(\Delta)[\mu_\Delta]$ contains a new set of reals.
If this were not the case, then $L(\Delta)[\mu_\Delta]$ would satisfy $\AD^+ + \DC + \neg\AD_\R+$``$\omega_1$ is $\ps(\R)$-supercompact'', which will be shown inconsistent in Section \ref{Club Filter}.
This means that in the first subcase of the second case, we also have $\kappa_\Delta = 0$.

Thus, the remaining subcase is when $\Delta$ is an $\AD_\R$ pointclass.
We place ourselves in $\Mod$ and argue as follows.
Let $(K, \Psi)$ be a mouse pair which is not in $\Delta$ and is minimal in the mouse order with respect to such mouse pairs.
In Section \ref{Description of a Generator}, we show that $K$ has a largest cardinal and is coded by its powerset.
Let $\M_\infty(K,\Psi)$ denote the direct limit of all non-dropping iterates of $(K,\Psi)$ and let $\kappa$ be the height of this direct limit.
We will show that $\kappa_\Delta \leq \kappa$ by establishing that $L(\Delta, \kappa^\omega)[\mu_\kappa]$ contains a set of reals which is not in $\Delta$.

This step is broken into two parts.
Let $S$ be the set of all $\omega$-sequences of elements of $\M_\infty(K,\Psi)$.
In Section \ref{Realizability Strategy}, we show that $L(\Delta, S)$ contains the non-dropping part of $\Psi$, which is essentially a set of reals not in $\Delta$.
This further reduces our goal to showing that $S \in L(\Delta, \kappa^\omega)[\mu_\kappa]$ and this is done in Section \ref{Adding Sequences and Measures}.
In Section \ref{Synthesis}, we put everything together and formally state and prove the main theorem.

\subsection{History}

Woodin has constructed\footnote{in the sense of that we described in Subsection \ref{Statement of the Result}} in an unpublished\footnote{but see \cite[Section 3]{steel2008long}} work the minimal model of $\AD_\R+\DC$.
Steel extended this work in \cite[Theorem 3.1]{steel2008long} and reached the minimal model of $\AD_\R+\DC+$``$\omega_1$ is $\ps(\R)$-supercompact''.
This latter model does not satisfy $V=L(\ps(\R))$, but its powerset of reals is presumably strictly larger than the powerset of reals in the former model.\footnote{We could not locate a proof of this claim in the literature.}
The construction that Steel proposes roughly corresponds to the first two oracles of our construction: if $\Delta$ is the current determinacy pointclass, then if $\cof(\Theta^{L(\Delta)})=\omega$, the next pointclass is
$$\ps(\R) \cap L(\Delta^\omega),$$
otherwise, the next pointclass is
$$\ps(\R) \cap L(\Delta)[\mu_\Delta],$$
but only in the case that it is strictly bigger then $\Delta$.
When this step fails to produce a bigger pointclass or produces a non-determined set, the construction terminates.
In the cited paper, Steel argues that the construction does not terminate too soon.
The pointclass that Steel's construction reaches is strictly smaller than the powerset of reals of the minimal model of $\AD_\R+$``$\Theta$ is regular''.\footnote{
To see this, work in $\AD_\R+$``$\Theta$ is regular''.
It turns out that the pointclass, letting $\Delta:=\Delta_{\theta_{\omega_1}}$, the model $L(\Delta)[\mu_\Delta]$ satisfies $\AD_\R+\DC+$``$\omega_1$ is $\ps(\R)$-supercompact'' and its powerset of reals is just $\Delta$.
Once we verify the latter, then the former follows by the results of \cite{solovay1978independence}.
The verification is based on the observation that every set of reals in $L(\Delta)[\mu_\Delta]$ is ordinal definable from a member of $\Delta$ and thus, in $\Delta$.
}

\subsection{Future Prospects}

Our construction produces the powersets of the reals of all $\AD_\R$ models below a strategy for a hod mouse with a measurable limit of Woodin cardinals.
We believe that this is not the ultimate extent of the method and that the construction reaches even larger $\AD_\R$ models.
However, our construction can reach at most the Chang-plus model
$$\mathsf{CM}^+:=L(\Ord^\omega)[ (\mu_{\alpha^\omega} : \alpha\in\Ord) ].$$
Namely, if there exists a proper class of Woodin limit of Woodin cardinals, the Chang-plus model satisfies $\AD_\R+\DC+$``$\omega_1$ is supercompact''.\footnote{
It follows from \cite[Theorem 1.3]{gappo2023chang} that it satisfies $\AD^++$``$\omega_1$ is supercompact'', while the rest follows by results of \cite{ikegami2021supercompactness}.
}
Moreover, the proof of the supercompactness of $\omega_1$ (due to Woodin) shows that it is witnessed by the restrictions of the club filters $\mu_X$ (for all $X$).
Letting $\Delta$ be the powerset of reals in the Chang-plus model, this is enough to conclude that for all $\kappa$,
$$\ps(\R) \cap L(\Delta^\omega, \kappa^\omega)[\mu_\Delta, \mu_\kappa] = \Delta.$$
In other words, if our construction can reach $\Delta$, it cannot go beyond it.
This means that a new kind of oracle is required at this stage.
Let us observe the following lower bound on the consistency of the determinacy beyond the Chang-plus model.

\begin{corollary}
    Suppose that $\mathsf{CM}^+$ satisfies $\AD^++$``$\omega_1$ is supercompact'' and that there exists an inner model $\Mod$ such that
    $$\ps(\R) \cap \mathsf{CM}^+ \subsetneq \Mod \models \AD^++\AD_\R.$$
    Then in $\Mod$, there exists a least branch hod pair with a measurable limit of Woodin cardinals.
\end{corollary}
\begin{proof}
    If there is such pair in $\mathsf{CM}^+$, we are done.
    Otherwise, the Main Theorem implies that there exists $\alpha < \sh$ such that
    $$\ps(\R) \cap \mathsf{CM}^+ = \De_\alpha.$$
    By the observations given above, our construction cannot continue further, so $\sh = \alpha+1$.
    Applying the Main Theorem this time to $\ps(\R) \cap \Mod$ shows that there must exist in $\Mod$ a least branch hod pair with a measurable limit of Woodin cardinals.
\end{proof}

Another missing point of our construction are the powersets of the reals of determinacy models satisfying $\AD^++\neg\AD_\R$.
Looking at our argument, one thing that we lack is the HOD analysis of these models in terms of least branch hod mice.
However, the conceptual difficulty lies in the usage of club filters.
Under $\AD_\R$, they are ultrafilters on any set which is a surjective image of the reals, but this may fail in the absence of $\AD_\R$.

Let us also say a few words on possible applications.
An important question that we would like to address is, working in some extension of $\ZFC$,\footnote{
for example, $\ZFC+\PFA$
} show that there are models of $\ZFC$ with a given large cardinal.\footnote{
In the $\PFA$ case, that can be a Woodin limit of Woodin.
Of course, this would require that our method is extended further.
}
One approach is to first exhibit an inner model containing all reals and satisfying $\AD^+$ and then show that its $\HOD$, or a rank initial segment thereof, has the target large cardinal.
Our method shows what this model must be, so what remains to show is that the construction lasts long enough as to ensure that $\HOD$ of the corresponding model is complex enough.\footnote{
In this scenario, the study of $\HOD$ would probably require the HOD analysis, which would require showing $\HPC$ in bigger and bigger determinacy models.
}

\subsection{Some Technical Remarks}

For the sake of precision, let us mentioned that by an \textit{inner model} we mean a transitive, proper class model satisfying $\ZF$.
We will be mostly interested in inner models containing all reals and satisfying $\AD^+$, which we might sometimes call \textit{determinacy models}.\footnote{
The intended reader of this paper is assumed to be familiar with the basic theory of $\AD^+$ and $\AD_\R$, but as a general reference we recommend \cite{larson2023extensions}.
}
Every such determinacy model contains $L(\R)$.
If $M \subseteq N$ are determinacy models, then $\ps(\R)^M \subsetneq \ps(\R)^N$ if and only if $\Theta^M < \Theta^N$.
Also, if $N$ is a deterimancy model and $M\subseteq N$ is an inner model containing all reals, then $M$ is a deterimancy model as well.\footnote{
This follows from \cite[Theorem 8.22]{larson2023extensions}.
The way this theorem is stated, it would require that $M$ is a definable subclass of $N$, but this issue can be avoid by assuming that both $M$ and $N$ satisfy $V=L(\ps(\R))$.
The reason why we can add this assumption is that an inner model $P$ satisfies $\AD^+$ if and only if $L(\ps(\R)^P)$ does so.
}

The construction that is the central topic of this paper concerns sufficiently closed pointclasses.
For our purposes, a pointclass is sufficiently closed if it is the collection of all sets of reals of some inner model containing all reals.
If $\Delta_0 \subseteq \Delta_1$ are closed pointclasses and $L(\Delta_1)$ satisfies $\AD^+$, then $L(\Delta_0)$ satisfies $\AD^+$ as well and it holds that $\Delta_0 \subsetneq \Delta_1$ if and only if
$$\Theta^{L(\Delta_0)} < \Theta^{L(\Delta_1)}.$$

\begin{notation}
    A subset $\Delta$ of $\ps(\R)$ is said to be a \intro{closed pointclass} iff 
    $$\ps(\R)\cap L(\R,\Delta)=\Delta.$$
    If $\Delta$ is a closed pointclass, then we denote by \intro{$\Theta_\Delta$} the ordinal $\Theta^{L(\Delta)}$.
    \qed
\end{notation}

Let us now fix some relatively standard notation in the $\AD$ context that we will use through the paper.

\begin{notation}
    Let us assume $\AD$.
    \begin{parts}
        \item We denote by \intro{$\Omega$} and \intro{$\theta_\alpha$} (for $\alpha \leq \Omega$) the length and the members of the Solovay sequence\footnote{
        Cf. the beginning of \cite[Section 6.3]{larson2023extensions}.
        }, respectively.\footnote{
        So in particular, $\theta_\Omega = \Theta$.
        }

        \item For $A\subseteq\R$, we denote by \intro{$\wa(A)$} the Wadge rank of $A$.
        
        \item For $\alpha \leq \Theta$, \intro{$\Delta_\alpha$} denotes the set of all $A\subseteq\R$ satisfying $\wa(A)<\alpha$.
        \qed
    \end{parts}
\end{notation}

For the fine structure, hod mice, and HOD analysis, our main sources are \cite{steel2023comparison} and \cite{steel2023mouse}.
We will try to follow their notation as much as possible.
One point where we diverge is the notation for the pull-back strategy.

\begin{notation}
    Let $M$ be a premouse, let $(N,\Tau)$ be a mouse pair, and let $\pi: M\to N$ be nearly elementary.
    Then we denote by \intro{$\pi^*\Tau$} the pull-back strategy of $\Tau$ along $\pi$.
    \qed
\end{notation}

Another detail that needs to be adjusted concerns cutpoints.

\begin{definition}
    Let $M$ be a premouse and let $\eta\leq o(M)$.
    Then we say that
    \begin{parts}
        \item $\eta$ is a \intro{cutpoint} of $M$ iff for all $\xi<\eta$, $o^M(\xi)\leq\eta$,
    
        \item $\eta$ is a \intro{strong cutpoint} of $M$ iff for all $\xi<\eta$, $o^M(\xi)<\eta$,
    
        \item $\eta$ is a \intro{clean cutpoint} of $M$ iff for all $M$-extenders $E$, if $\len(E)>o(M)$, then $\crit(E)>o(M)$.
        \qed
    \end{parts}
\end{definition}

Cutpoints and strong cutpoints appear in \cite[Definition 2.16]{steel2023mouse}, while the clean cutpoints are new.
The following relation is the reason why we introduce them.

\begin{definition}
    Let $M$ and $N$ be premice.
    Then \intro{$M\unlhd^*N$} iff $M\unlhd N$ and $o(M)$ is a clean cutpoint of $N$.
    \qed
\end{definition}

This relation corresponds to the strong initial segment relation from \cite[Definition 4.1]{steel2023mouse}, but the way the relation is defined there excludes the possibility that $M\unlhd^* N$ whenever $M$ is extender-active.
We believe this to be a typo and that our formulation captures the point.

All the notation that we introduce throughout the paper is compiled in the index at the end.

\newpage

\section{HOD Mice and HOD Analysis}
\label{HOD Mice and HOD Analysis}

We try to follow \cite{steel2023comparison} regarding the fine structure as much as possible.
To avoid cumbersome wording, we will use the following convention.

\begin{notation}
    When we say ``\intro{premouse}'' and ``\intro{mouse pair}'', we mean the one of the least branch hod type.
    The scope of a mouse pair is assumed to be $H_{\omega_1}$, unless stated otherwise.
    \qed
\end{notation}

\begin{notation}
    Assume $\AD^+$.
    Let $(M,\Sigma)$ be a mouse pair.
    We denote by \intro{$(M,\Sigma) \para \xi$} the mouse pair
    $$(M\para\xi,\Sigma_{M\para\xi}),$$
    whenever $\xi\leq\hat o(M)$.
    The pairs \intro{$(M,\Sigma)\vert\xi$} and \intro{$(M,\Sigma)\vert(\xi,n)$} are defined similarly (for $(\xi,n)\leq l(M)$).
    \qed
\end{notation}

We consider that every notion defined for premice is inherited by mouse pairs via their premouse component, unless explicitly stated otherwise.
For example, a mouse pair $(M,\Sigma)$ is strongly stable iff $M$ is strongly stable, where the definition of strong stability of premice reads as follows.

\begin{notation}
    A premouse $M$ is \intro{strongly stable} iff the $r\Sigma_{k(M)}$-cofinality of $\rho_{k(M)}(M)$ is not measurable by the $M$-sequence.
\end{notation}

The comparison theorem that we work with is proved for strongly stable mouse pairs, cf. \cite[Theorem 9.3.6]{steel2023comparison}.
The shows the importance of strongly stable mouse pairs.
In particular, given a mouse pair $(M,\Sigma)$ which is \underline{not} strongly stable, we will often have to replace it with one that is so.
The construction that we most often use in this paper comes from \cite[Page 111, last paragraph]{steel2023comparison} and \cite[Page 43, penultimate paragraph]{steel2023mouse} and it goes as follows.
Since $M$ is not strongly stable, we necessarily have that $k:=k(M)>0$.
We denote by $\Bar{M}$ the premouse
$$\mathfrak{C}(M^-)$$
and by $\Bar{\Sigma}$ the pullback of $\Sigma$ by the anti-core embedding.
Then $(\Bar{M},\Bar{\Sigma})$ is essentially a mouse pair, except that $\Bar{M}$ is may fail to be pfs.
The only thing that can fail is in fact that its $k^\mathrm{th}$ projectum is measurable.
In other words, pair $(\Bar{M},\Bar{\Sigma})$ is either equal to $(M,\Sigma)$ or $(M,\Sigma)^-$ is the ultrapower of $(\Bar{M},\Bar{\Sigma})^-$ by an order zero measure on $\rho_k(\Bar{M})$.
Since $M$ is assumed not to be strongly stable, there is an order zero measure $D$ on the $r\Sigma_k$-cofinality of $\rho_k(M)$.
Let $(M^*,\Sigma^*)$ be the ultrapower of $(\Bar{M},\Bar{\Sigma})$ via $D$.
Then $(M^*,\Sigma^*)$ is a strongly stable mouse pair.

One drawback of the above construction is that if $\Bar{M}$ is not equal to $M$, then $(M,\Sigma)^-$ is an ultrapower of $(\Bar{M},\Bar{\Sigma})^-$, but it is not clear that we can reconstruct the full strategy $\Sigma$ from $\Bar{\Sigma}$.
This prevents us from showing that when passing from $\Sigma$ to $\Sigma^*$ we do not drop in Wadge rank.
However, there is another way to make a mouse pair strongly stable and it is based on \cite[Section 3]{steel2023condensation}.\footnote{
We thank John Steel for drawing our attention to this.
}
Since we are assuming that $M$ is not strongly stable, we can take the order zero measure $D$ on $\eta_k^M$ and produce the ultrapower $\Ult(M,D)$.
An issue is that the premouse $\Ult(M,D)$ might not be $k$-sound, i.e. it might be a pfs premouse of type 2.
If that happens, we have to make sense of iterations of such non-sound structures in order to talk about the mouse pair $(\Ult(M,D),\Sigma_{(D)})$.
This is done in \cite[Lemma 3.49]{steel2023condensation}.
What we care about is here is the $k^\mathrm{th}$ core $M^*$ of $\Ult(M,D)$, together with the pullback strategy $\Sigma^*$ of $\Sigma_{(D)}$ along the anti-core map.
This is a strongly stable mouse pair (of type 1).
The cited lemma shows that $\wa(\Sigma_{(D)})\geq\wa(\Sigma^*)$, while the reverse inequality follows is witnessed simply by the anti-core map.
Since strategies $\Sigma$ and $\Sigma_{(D)}$ are an $r\Sigma_k$-ultrapower away from each other, they also have the same Wadge-rank, finally ensuring the conclusion that $\wa(\Sigma^*)=\wa(\Sigma)$.

One way in which this construction turns out to be useful is with respect to $\HPC$.
We invite the reader to consult \cite[Definition 1.7.1]{steel2023mouse} and the paragraph that follows if they are not familiar with this principle.
As stated there, assuming $\AD^+$, principle $\HPC$ is equivalent to the strategies of mouse pairs being Wadge-cofinal in Suslin co-Suslin sets of reals.
What the above construction brings to the table is that we can make sure that these mouse pairs are strongly stable, which gives the following lemma.

\begin{lemma}
    Suppose that $\AD^+$ holds.
    Then $\HPC$ holds if and only if for all Suslin co-Suslin sets of reals $A$, there exists a strongly stable mouse pair $(M,\Sigma)$ such that $\wa(\Sigma)\geq\wa(A)$.\footnote{
    This was stated without an explanation in the paragraph following \cite[Definition 0.6]{steel2023mouse}.
    }
    \qed
\end{lemma}

The principle $\HPC$ is conjectured to follow from $\AD^+$ and the assumption that there are no mouse pairs with long extenders.\footnote{
Cf. \cite[Conjecture 1.7.6]{steel2023comparison}.
}
If we (significantly) strengthen this additional requirement, then $\HPC$ is known to hold.
The anti-large-cardinal assumption that will suffice for our purposes reads as follows.

\begin{notation}
    The statement \intro{$\mathsf{NMLW}$} (``no measurable limit of Woodins'') asserts that ``there does \underline{not} exist a mouse pair $(M,\Sigma)$ such that $M$ has a measurable cardinal which is a limit of Woodin cardinals''.
    \qed
\end{notation}

The important theorem which will often be implicitly used is the following.

\begin{theorem}\label{HPC}
    $\mathsf{ZF}+\AD^++\mathsf{NMLW}\vdash\HPC$
\end{theorem}
\begin{proof}
    Working in $\ZF+\AD^++\mathsf{NMLW}$, we know by \cite[Theorem 10.2.1]{sargsyan2023largest} that the \textit{Mouse Capturing} holds.
    It was shown in \cite{steel2008derived} that the \textit{Mouse Capturing} is equivalent to the principle $\mathsf{LEC}$ defined in \cite[Defition 1.7.1]{steel2023comparison}.
    By \cite[Defition 10.4.3]{steel2023comparison}, $\mathsf{LEC}$ implies $\HPC$.
\end{proof}

The importance of $\HPC$ lies in the fact that it allows \textit{HOD analysis}, which will be the main technical tool in the paper.
Our starting point for that method is \cite{steel2023mouse} and we will summarize here the aspects of the paper that will be important for us.

\begin{definition}\label{generator}
    Assume $\AD^+$.
    Let $\Delta$ be a selfdual pointclass and let $(M,\Sigma)$ be a mouse pair.
    Then $(M,\Sigma)$ is a \intro{generator} for $\Delta$ iff
    \begin{parts}
        \item $(M,\Sigma)$ is strongly stable,
        
        \item $\Sigma$ is not in $\Delta$,

        \item for all mouse pairs $(N,\Tau)<^*(P,\Sigma)$, we have that $\Tau$ is in $\Delta$.
        \qed
    \end{parts}
\end{definition}

We impose the strong stability on the generators from the beginning.\footnote{
In \cite[Definition 3.3]{steel2023mouse}, the author did not require the strong stability, which we believe to be a non-essential omission.
}
The reason for this is that strongly stable mouse pairs are pre-well-ordered by the mouse order.\footnote{
Cf. \cite[Corollary 9.3.7.b]{steel2023comparison}.
}
Hence, a generator for $\Delta$ is the least mouse pair outside of $\Delta$ and all the generators for $\Delta$ are mouse-equivalent.
When generating pointclasses, essentially only Suslin co-Suslin sets matter.

\begin{lemma}
    Let $\Delta$ be a closed pointclass and let $\Delta'$ be the pointclass of Suslin co-Suslin sets in $L(\Delta)$.
    Then $(M,\Sigma)$ is a generator for $\Delta$ if and only if it is a generator for $\Delta'$.
\end{lemma}
\begin{proof}
    This follows from the fact that the iteration strategies are Suslin co-Suslin.
\end{proof}

The mouse limits of generators for non-limit Solovay pointclasses are known full in the following sense.

\begin{lemma}\label{174}
    Suppose that
    \begin{assume}
        \item $V=L(\ps(\R))\models\AD^+ + \AD_\R+\HPC$,
        
        \item $\alpha<\Omega$ is \underline{not} a limit,

        \item $(M,\Sigma)$ is a generator for $\Delta_{\theta_\alpha}$,

        \item $(N,\Tau)$ is a strongly stable mouse pair,

        \item $(M,\Sigma)\leq^* (N,\Tau)$.
    \end{assume}
    Then $\M_\infty(M,\Sigma)\unlhd^*\M_\infty(N,\Tau)$.
\end{lemma}
\begin{proof}
    This is \cite[Corollary 4.5]{steel2023mouse}.
\end{proof}

This fullness allows us to stack all such mouse limits.

\begin{definition}\label{227}
    Assume $V=L(\ps(\R))\models\AD^+ + \AD_\R+\HPC$.
    We define \intro{$\hh$} to be the least branch premouse of the least possible length satisfying that for all $\xi<\Omega$, for all $(M,\Sigma)$ which are generators for $\Delta_{\theta_{\xi+1}}$, we have that $\M_\infty(M,\Sigma)\unlhd^*\hh$.
    \qed
\end{definition}

In turns out that this stack is exactly the $\HOD$ of the determinacy model.

\begin{theorem}\label{what is hod}
    Assume $V=L(\ps(\R))\models\AD^+ + \AD_\R+\HPC$.
    The domain of $\hh$ is equal to $\HOD\para\Theta$.
\end{theorem}
\begin{proof}
    This is \cite[Theorem 4.7]{steel2023mouse}.
\end{proof}

\begin{notation}
    Let $\Gamma$ be a closed pointclass such that $L(\Gamma)\models\AD^+ + \AD_\R+\HPC$.
    Then we denote by \intro{$\hh_\Gamma$} the premouse $\hh^{L(\Gamma)}$.
    \qed
\end{notation}

We will also need some additional properties of generators of the non-limit Solovay pointclasses.
Before we state them, let us recall some notation.

\begin{definition}
    Suppose that $M$ is a premouse.
    Then we define the following.
    \begin{parts}
        \item The ordinal \intro{$\tau^M$} is the supremum of all $\kappa^{+M}$ where $\kappa<\rho^-(M)$ and there exists $\eta\in (\kappa^{+M},o(M)]$ such that $\kappa=\crit(E^M_\eta)$.

        \item If $\tau^M$ is not a strong cutpoint of $M$, then the ordinal \intro{$\beta^M$} the least ordinal $\xi<o(M)$ satisfying that $o(\xi)^M\geq\tau^M$.\footnote{
        We remind the reader that $o(\xi)^M$ denotes the ordinal $\sup\{\eta\leq o(M) : \crit(E^M_\eta)=\xi\}$.
        }
        \qed
    \end{parts}
\end{definition}

\begin{definition}
    Let $(M,\Sigma)$ be a strongly stable mouse pair.
    Then we define the following.
    \begin{parts}
        \item The ordinal \intro{$\tau_\infty(M,\Sigma)$} is defined as $\pi_{(M,\Sigma),\infty}(\tau^M)$.

        \item If $\tau^M$ is not a strong cutpoint of $M$, then the ordinal \intro{$\beta_\infty(M,\Sigma)$} is defined as $\pi_{(M,\Sigma),\infty}(\beta^M)$.
        \qed
    \end{parts}
\end{definition}

The properties that we want now read as follows.

\begin{theorem}\label{generator Solovay successor}
    Suppose that
    \begin{assume}
        \item $V=L(\ps(\R))\models\AD^+ + \AD_\R+\HPC$,
        
        \item $\alpha<\Omega$ is \underline{not} a limit,

        \item $(M,\Sigma)$ is a generator for $\Delta_{\theta_{\alpha}}$,

        \item $\Delta$ is the pointclass of Suslin co-Suslin sets of $L(\Delta_{\theta_{\alpha}})$.
    \end{assume}
    Then it holds that
    \begin{parts}
        \item $o(M)=\tau^M$ and $k(M)=0$,

        \item\label{193} $M$ is passive,


        \item the domain of $\M_\infty(M,\Sigma)$ is $\HOD\para\theta_\alpha$,

        \item $\beta_\infty(M,\Sigma)=o(\Delta)$ and $\tau_\infty(M,\Sigma)=\theta_\alpha$.
    \end{parts}
\end{theorem}
\begin{proof}
    These facts are listed in the first paragraph of \cite[Page 38]{steel2023mouse}, except \ref{193}.
    To see \ref{193}, fix a generator $(N,\Tau)$ for $\Delta_{\theta_{\alpha+1}}$.
    We have that $o(\M_\infty(N,\Tau))=\theta_{\alpha+1}$.
    By Lemma \ref{174}, we get that $\M_\infty(M,\Sigma)\lhd\M_\infty(N,\Tau)$.
    Since $\theta_\alpha=o(\M_\infty(M,\Sigma))$ is a cardinal in $V$, it is so in $\M_\infty(N,\Tau)$ as well.
    This means that in $\M_\infty(N,\Tau)$, neither a branch nor an extender are indexed at $o(\M_\infty(M,\Sigma))$, which shows that $\M_\infty(M,\Sigma)$ is passive.
\end{proof}

\newpage

\section{Club Filter}
\label{Club Filter}

The goal of this section is to describe the notion of club filter in the choiceless context.
This notion will generalize the usual notion of club filter under the Axiom of Choice.
The definition that we will work with essentially comes from \cite{solovay1978independence}.

\begin{declaration}
    Suppose that
    \begin{assume}
        \item $\ZF$ holds,

        \item $S$ is an infinite set,

        \item $\DC_S$ holds.
        \qed
    \end{assume}
\end{declaration}

Under these assumptions we know that the set $[S]^\omega$ of all countable subsets of $S$ is non-empty.
We want to talk about filters on $[S]^\omega$, one instance of which is the club filter.
Below are some properties which are usually associated with the club filter.


\begin{definition}
    Let $\ff$ be a filter on $[S]^\omega$.
    \begin{parts}
        \item $\ff$ is said to be \intro{fine} iff for all $\tau\in [S]^\omega$,
        $$\{\sigma\in [S]^\omega : \tau\subseteq\sigma\}\in\ff.$$

        \item $\ff$ is said to be \intro{normal} iff for all $A$ which are $\ff$ positive, for all
        $$F:A\to [S]^\omega$$
        satisfying that
        $$\forall\sigma\in A, \emptyset\subset F(\sigma)\subseteq\sigma,$$
        there exists $a\in S$ such that
        $$\{\sigma\in A : a\in F(\sigma)\}$$
        is $\ff$-positive.

        \item $\ff$ is said to have the \intro{diagonal intersection property} iff for all families $(A_a : a\in S)$ of sets in $\ff$, the set
        $$\mathop{\Delta}_{a\in S}A_a:=\{\sigma\in [S]^\omega : \forall a\in\sigma, \sigma\in A_a\}$$
        belongs to $\ff$.
        \qed
    \end{parts}
\end{definition}

The notion of normality is adjusted to the choiceless context.
This adjustment is appropriate, since we have the following standard equivalence.

\begin{proposition}
    Let $\ff$ be a filter on $[S]^\omega$.
    Then $\ff$ is normal iff $\ff$ has the diagonal intersection property.
\end{proposition}
\begin{proof}
    By the proof of \cite[Lemma 4.5]{solovay1978independence}.
\end{proof}

We also have that normal fine filters are $\omega_1$-complete.

\begin{lemma}\label{389}
    Suppose that $\ff$ is a normal, fine filter on $[S]^\omega$.
    Then $\ff$ is $\omega_1$-complete.
\end{lemma}
\begin{proof}
    Let $(A_n : n<\omega)\in\ff^\omega$ be arbitrary.
    There exists an injection $(a_n : n<\omega)\in S^\omega$.
    We define $(B_a : a\in S)\in\ff^S$ as follows: for $a\in S$,
    $$B_a:=\begin{cases}
        A_n, & n<\omega, a=a_n\\
        S, & \mbox{otherwise.}
    \end{cases}$$
    Since $\ff$ is fine, we have that
    $$C:=\{\sigma\in [S]^\omega : \{a_n\}_{n<\omega}\subseteq\sigma\}\in\ff.$$
    Note that
    $$\bigcap_{n<\omega}A_n\cap C=\mathop{\Delta}_{a\in S}B_a\cap C.$$
    Since $\ff$ is normal, it follows that $\bigcap_{n<\omega}A_n\in\ff$.
\end{proof}

Following the terminology \cite{vaananen2011models}, we introduce the so-called \textit{club game}.
This game is what will allow us to characterize the club filter.

\begin{definition}
    Let $A\subseteq [S]^\omega$.
    The game \intro{$\Game^\mathrm{club}(S,A)$} is a length $\omega$ two player game of the form
    \begin{center}
        \begin{tabular}{c|ccccc}
            I  & $a_0$ &       & $a_1$ & & $\cdots$\\
            \hline
            II &       & $b_0$ & & $b_1$ & $\cdots$
        \end{tabular}
    \end{center}
    where for all $n<\omega$, $a_n,b_n\in S$.
    Player II wins iff $\{a_n,b_n : n<\omega\}\in A$.
    \qed
\end{definition}

We will introduce now what is supposed to be the club filter, but we will avoid using the word ``club'' unless the Axiom of Choice holds.
The reason for this is that different characterizations of being a club are not necessarily equivalent in the choiceless context.
In informal discussions, we might use the term ``quasi-club filter''.

\begin{definition}\label{199}
    We define \intro{$\mu_S$} to consist of all $A\subseteq [S]^\omega$ for which Player II has a winning quasi-strategy in $\Game^\mathrm{club}(S,A)$.
\end{definition}

The set $\mu_S$ is a filter and it will have some of the usual properties of the club filter.

\begin{proposition}
    $\mu_S$ is a fine filter.
\end{proposition}
\begin{proof}
    To see that $\mu_S$ is a filter, the only non-trivial part is to verify the closure for intersections.
    This is argued similarly to the countable completeness of \cite[Lemma 4.1]{solovay1978independence}, but we work with quasi-strategies instead of strategies.
    Our case is simpler insomuch that we are only considering intersections of two sets.
    Let us explain why $\mu_S$ is fine.
    Let $\sigma\in [S]^\omega$ and let
    $$A:=\{\sigma'\in [S]^\omega : \sigma\subseteq\sigma'\}.$$
    We need to give a winning quasi-strategy for Player II in $\Game^\mathrm{club}(S,A)$.
    If we fix an enumeration $(a_n : n<\omega)$ of $\sigma$, then Player II can win $\Game^\mathrm{club}(S,A)$ by playing $a_n$'s, regardless of the moves of Player I.
\end{proof}

The normality seems to require more than just the assumptions that we are working under.
One way to ensure it is to assume more choice.

\begin{proposition}
    Suppose that $\DC_{\ps(S^{<\omega})}$ holds.
    Then filter $\mu_S$ is normal.
\end{proposition}
\begin{proof}
    Suppose that $(A_a : a\in S)$ is a family of sets in $\mu_S$. 
    Implicit in the proof of \cite[Lemma 4.6]{solovay1978independence} is a quasi-strategy $\tau$ for Player II in the game in the game
    $$\Game^\mathrm{club}(S,\mathop{\Delta}_{a\in S}A_a).$$
    Note that we do not pick a family $(\tau_a : a\in S)$ of winning quasi-strategy for Player II in the games
    $$(\Game^\mathrm{club}(S,A_a) : a\in A)$$
    as it is done there, but we play a tree of auxiliary games according to all possible such quasi-strategies.
    In the end, if $f$ is a play according to $\tau$, we can use $\DC_{\ps(S^{<\omega})}$ to be pick a sequences of $(\tau_a : a\in\ran(f))$ of winning quasi-strategies for Player II in the games
    $$(\Game^\mathrm{club}(S,A_a) : a\in \ran(f))$$
    which ensure that for all $a\in\ran(f)$, $\ran(f)\in A_a$.
\end{proof}

We will soon observe another natural way to obtain the normality.
In the choiceless context that we care about the most, the quasi-club filter will be both ultra and normal.
We call such filters \textit{supercompactness measures}.

\begin{definition}
    A \intro{supercompactness measure} on $[S]^\omega$ is a normal, fine ultrafilter on $[S]^\omega$.
    \qed
\end{definition}

It follows from Lemma \ref{389} that every supercompactness measure on $[S]^\omega$ is $\omega_1$-complete.
We will now show that every supercompactness measure must extend the quasi-club filter.

\begin{proposition}\label{cf-149}
    Suppose that $\mu$ is a supercompactness measure on $[S]^\omega$.
    Then $\mu_S\subseteq\mu$.
\end{proposition}
\begin{proof}
    This follows from the proof of \cite[Theorem 1.2]{woodin2021ad}.
\end{proof}

One immediate consequence is that the existence of a supercompactness measure implies the normality of the quasi-club filter.
Another consequence is that if the quasi-club filter itself is a supercompactness measure, then it is the only supercompactness measure.
This situation naturally occurs in the determinacy context (see Proposition \ref{485} below).

\begin{lemma}\label{cf-290}
    Suppose that $\AD_\R$ holds and that $\R$ surjects onto $S$.
    Then in every game on $S$, one of the players has a winning quasi-strategy.
\end{lemma}
\begin{prooff}
    \item Let $A\subseteq S$ be arbitrary.
    We want to show that one of the players in the game $\Game(S,A)$ on $S$ with $A$ as the payoff set for Player II has a winning quasi-strategy.

    \item Let $\pi:\R\twoheadrightarrow S$ be a surjection and let
    $$B:=\{y\in\R^\omega : \pi\circ y\in A\}.$$
    By $\AD_\R$, one of the players in the game on $\R$ with $B$ as the payoff set for Player II has a winning strategy.
    We can assume without loss of generality that it is Player I who has a winning strategy $\tau$.

    \item Regarding $\tau$ as a tree, let $\sigma:=\{\pi\circ y : y\in\tau\}$.
    Then $\sigma$ is a quasi-strategy for Player I in $\Game(S,A)$.
    We want to show that it is winning.

    \item Let $f$ be a play according to $\sigma$.
    By $\DC_\R$, there exists $g$ according to $\tau$ such that $f=\pi\circ g$.

    \item Since $g$ is according to $\tau$, we have that $g\in B^\complement$.
    This implies $f\in A^\complement$, as required.
\end{prooff}

\begin{proposition}\label{485}
    Suppose that $\AD_\R$ holds and that $\R$ surjects onto $S$.
    Then $\mu_S$ is a supercompactness measure.
    Furthermore, it is the only supercompactness measure on $[S]^\omega$.\footnote{
    \cite[Theorem 1.5]{woodin2021ad} proves that there exists a unique supercompactness measure on $[S]^\omega$, but we need a more precise statement.}
\end{proposition}
\begin{proof}
    By Lemma \ref{cf-290}, either Player I or Player II has a winning quasi-strategy in $\Game^\mathrm{club}(S,A)$.
    If Player II has a winning quasi-strategy, then $A\in\mu_S$ by the definition.
    On the other hand, if Player I has a winning quasi-strategy in $\Game^\mathrm{club}(S,A)$, then we can easily show that Player II has a winning quasi-strategy in $\Game^\mathrm{club}(S,A^\complement)$, which means that $A^\complement\in\mu_S$.
    Thus, $\mu_S$ is a fine ultrafilter on $[S]^\omega$.
    To see that it is also normal, we can run the proof of \cite[Lemma 4.6]{solovay1978independence}.
    On the other hand, the uniqueness follows from \ref{cf-149}.
\end{proof}

In the case that $S$ is a mouse limit, there is an important collection of quasi-clubs that we will have to be aware of in the proof of Theorem \ref{1882}.

\begin{proposition}\label{535}
    Suppose that
    \begin{assume}
        \item $\AD^+$ holds,

        \item $(M,\Sigma)$ is a mouse pair,

        \item $S:=\M_\infty(M,\Sigma)$.\footnote{
        Note that in this case $\DC_S$ holds because $S$ is wellordered, so adding this assumption is not necessary.
        }
    \end{assume}
    Then the set $\{\ran (\pi_{(N,\Tau),\infty}) : (N,\Tau)\mbox{ is a non-dropping iterate of }(M,\Sigma)\}$ belongs to $\mu_S$.
\end{proposition}
\begin{proof}
    Let us denote by $A$ the said set.
    A winning quasi-strategy for Player II in $\Game^\mathrm{club}(S,A)$ consists in the following.
    Let $(M_{-1},\Sigma_{-1}):=(M,\Sigma)$.
    As the game progresses and Player I plays some $a_n\in S$ (for $n<\omega$), Player II produces a non-dropping iterate $(M_n, \Sigma_n)$ of $(M_{n-1}, \Sigma_{n-1})$ so that $a_n\in \ran(\pi_{(M_n,\Sigma_n),\infty})$.
    Player II then makes sure to list each element of $\ran(\pi_{(M_n,\Sigma_n),\infty})$ in one of his future moves.
    The set produced in the end is equal to
    $$\bigcup_{n<\omega}\ran(\pi_{(M_n,\Sigma_n),\infty}),$$
    which is easily seen to belong to $A$.
    One thing to keep in mind is that we actually never have to pick an iterate $(M_n, \Sigma_n)$ (or its enumeration in order type $\omega$) since we only need to exhibit a quasi-strategy.
\end{proof}

Another fact that we will use in the proof of Theorem \ref{1882} is that quasi-club filters project onto each other.

\begin{lemma}\label{562}
    Suppose that
    \begin{assume}
        \item $\Bar{S}\subseteq S$ is infinite,\footnote{
        Observe that $\DC_S$ implies $\DC_{\Bar{S}}$ since every tree on $\Bar{S}$ is also a tree on $S$.
        }

        \item $A$ belongs to $\mu_S$.
    \end{assume}
    Then the set $\{\sigma\cap\Bar{S} : \sigma\in A\}$ belongs to $\mu_{\Bar{S}}$.
\end{lemma}
\begin{proof}
    Let $u$ be an arbitrary element of $\Bar{S}$ and let $A'$ consist of all $\sigma\in A$ satisfying that $u \in \sigma$.
    Since $\mu_S$ if fine, it follows that $A'\in \mu_S$.
    Let us consider an arbitrary winning quasi-strategy for Player II in $\Game^\mathrm{club}(S,A')$.
    We transform this quasi-strategy into a quasi-strategy for Player II in 
    $$\Game^\mathrm{club}(\Bar{S},\{\sigma\cap\Bar{S} : \sigma\in A'\})$$
    by asserting that Player II should follow it faithful as long as it prescribes an element of $\Bar{S}$ and otherwise, he should play $u$ instead.
    This new quasi-strategy is easily seen to be winning, thus witnessing that
    $$\{\sigma\cap\Bar{S} : \sigma\in A'\} \in \mu_{\Bar{S}}.$$
    The conclusion now readily follows from the fact that $A\supseteq A'$.
\end{proof}

Let us state for the record that under the Axiom of Choice, the quasi-club filter is simply the club filter.

\begin{proposition}
    Suppose that $\mathsf{AC}$ holds.
    Then $\mu_S$ is equal to the club filter on $[S]^\omega$.
\end{proposition}
\begin{proof}
    $\mathsf{AC}$ allows us to work with strategies instead of quasi-strategies.
    The conclusion then follows from \cite[Proposition 6.19]{vaananen2011models}.
\end{proof}

From this point on, we work towards proving the Main Theorem.
The following proposition has for a consequence that the construction described in the introduction is sufficiently absolute.

\begin{proposition}\label{509}
    Suppose that
    \begin{assume}
        \item $M$ is an inner model,

        \item $S\in M$,

        \item $S^\omega\subseteq M$.
    \end{assume}
    Then $\mu_S^M\subseteq\mu_S$.
    In particular, if $\mu_S^M$ is an ultrafilter in $M$, then $\mu_S^M=\mu_S\cap M$.
\end{proposition}
\begin{proof}
    Suppose that $A\in\mu_S^M$.
    Then there exists $\tau$ such that $M\models$``$\tau$ is a winning quasi-strategy for Player II in $\Game^\mathrm{club}(S,A)$''.
    Since $S^\omega\subseteq M$, we see that $\tau$ is a winning quasi-strategy for Player II in $\Game^\mathrm{club}(S,A)$ (in $V$).
    This shows that $A\in\mu_S$.
\end{proof}

In the introduction, we said that in the case that $\Delta$ is a deterimancy pointclass which is \underline{not} an $\AD_\R$ pointclass and for which $\cof(\Theta^{L(\Delta)})>\omega$, the model $L(\Delta)[\mu_\Delta]$ contains a set of reals which is not in $\Delta$.
This is the content of the corollary below, which states in the counter-positive.

\begin{proposition}\label{519}
    Suppose that $\AD^+$ and $\DC$ hold and that there exists a supercompactness measure on $[\ps(\R)]^\omega$.
    Then $\AD_\R$ holds.\footnote{
    \cite[Theorem 2.2.1]{trang2013generalized} shows that the conclusion does not follow if instead of assuming that $\omega_1$ is $\ps(\R)$-supercompact, we just assume that it is $\R$-supercompact.
    }
\end{proposition}
\begin{proof}
    Since we are assuming $\AD^+$, to get $\AD_\R$ it suffices to show that every Suslin set is co-Suslin.\footnote{
    Cf. \cite[Theorem 0.3]{larson2023extensions}.
    }
    Let $\lambda<\Theta$ and let $T$ be a tree on $\omega\times\lambda$.
    We want to show that $p[T]^\complement$ is Suslin.
    We can run the proof of \cite[Theorem 2]{ikegami2021supercompactness} to get that $T$ is weakly homogeneously Suslin.
    Since $\DC$ holds, we can construct a Martin-Solovay tree for $T$ and show that it projects onto $p[T]^\complement$.
    This tree witnesses that $p[T]^\complement$ is Suslin.

    We briefly explain why the cited proof goes through.
    Its examination shows that the authors are using that $\DC$ holds, that there exists a supercompactness measure on\footnote{Following the cited opus, $\mathsf{MEAS}^{\omega_1,\lambda}_n$ denotes the set of all countably complete measures on $\lambda^n$.}
    $$\left[\bigcup_{n<\omega}\left(\ps(\lambda^n)\cup\mathsf{MEAS}^{\omega_1,\lambda}_n\right)\right]^\omega,$$
    and that there exists a supercompactness measure on $[\ps(\lambda)]^\omega$.
    In our case, $\DC$ hods by the assumption.
    To get the required supercompactness measures, observe that $\lambda<\Theta$ means that $\R$ surjects to $\lambda$, so by Coding Lemma, $\R$ surjects onto $\ps(\lambda)$.
    This suffices to get the second supercompactness measure.
    To get the first one, we need to show that there is a surjection
    $$\ps(\R)\twoheadrightarrow\bigcup_{n<\omega}\left(\ps(\lambda^n)\cup\mathsf{MEAS}^{\omega_1,\lambda}_n\right).$$
    Since $\DC$ holds, it suffices to show that for each $n<\omega$, there exists a surjection
    $$\ps(\R)\twoheadrightarrow\ps(\lambda^n)\cup\mathsf{MEAS}^{\omega_1,\lambda}_n.$$
    Since there is a canonical bijection $\lambda\leftrightarrow\lambda^n$ for $n>1$, we may assume that $n=1$.
    In that case,
    $$\ps(\lambda^n)\cup\mathsf{MEAS}^{\omega_1,\lambda}_n\subseteq\ps(\lambda)\cup\ps(\ps(\lambda)),$$
    so the required surjection is obtained by the fact that $\R$ surjects onto $\ps(\lambda)$.
\end{proof}

\begin{corollary}\label{565}
    Suppose that $\AD^+$ holds, that $\mu_{\ps(\R)}$ is a supercompactness measure, and that $V=L(\ps(\R))[\mu_{\ps(\R)}]$.
    Then $V\models\AD_\R+\DC$.
\end{corollary}
\begin{proof}
    Note that $\omega_1$ is $\ps(\R)$-supercompact.
    The proof of \cite[Theorem 1]{ikegami2021supercompactness} then shows that $\DC_{\ps(\R)}$ holds.
    Since $V=L(\ps(\R))[\mu_{\ps(\R)}]$, we conclude that full $\DC$ holds.
    Proposition \ref{519} shows that $\AD_\R$ holds.
\end{proof}

\newpage

\section{Description of a Generator}
\label{Description of a Generator}

\begin{declaration}
    Suppose that
    \begin{assume}
        \item $V=L(\ps(\R))\models\AD^+ +\HPC$,

        \item there does \underline{not} exist a mouse pair with a superstrong cardinal,

        \item $\Gamma$ is a closed pointclass such that $L(\Gamma)\models\AD_\R$,

        \item $\omega<\cof(\Theta_\Gamma)\leq \Theta_\Gamma<\Theta$,

        \item $(K,\Psi)$ is a generator for $\Gamma$.
        \qed
    \end{assume}
\end{declaration}

The goal of this section is establish how the generator $K$ for $\Gamma$ looks like.
For the final description, we will have to strengthen the anti-large-cardinal assumption to the non-existence of a mouse pair with a measurable limit of Woodin cardinals.
In that case, $K$ has the largest cardinal, which is a limit of Woodin cardinals, and its fine structural degree is $0$.
To start arguing, we will require the following technical notion from \cite{steel2023mouse}.

\begin{definition}
    Suppose that $(M,\Sigma)\leq^*(N,\Tau)$ are strongly stable mouse pairs.
    Then $(\T,\uu)$ is a \intro{minimal comparison} of $(M,\Sigma)$ with $(N,\Tau)$ iff
    \begin{parts}
        \item $\T$ is a normal tree on $(M,\Sigma)$ with the last pair $(M',\Sigma')$,

        \item $\uu$ is a normal tree on $(N,\Tau)$ with the last pair $(N',\Tau')$,

        \item the main branch of $\T$ does not drop,

        \item $(M',\Sigma')\unlhd^*(N',\Tau')$,


        \item for all $\alpha<\len(\uu)$, $\crit(E^\uu_\alpha)\leq o(M')$.
        \qed
    \end{parts}
\end{definition}

\begin{lemma}
    Suppose that $(M,\Sigma)\leq^*(N,\Tau)$ are strongly stable mouse pairs.
    Then there exists a minimal comparison of $(M,\Sigma)$ with $(N,\Tau)$.
    \qed
\end{lemma}

We will also need a fullness lemma for generators of Solovay pointclasses.

\begin{lemma}\label{fullness of Solovay generators}
    Let $(M,\Sigma)$ be a generator for some $\Delta_{\theta_{\alpha+1}}^{L(\Gamma)}$, where $\alpha<\Omega^{L(\Gamma)}$.
    Then for all strong cutpoint cardinals $\kappa$ of $M$, letting $\Bar{M}:=M|\kappa^{+M}$, for all mouse pairs $(N,\Tau)\in L(\Gamma)$ satisfying that $(N,\Tau)\unrhd (\Bar{M},\Sigma_{\Bar{M}})$, it holds that that $\rho(N)>\beta^M$.
\end{lemma}
\begin{proof}
    This is \cite[Lemma 4.3]{steel2023mouse} relativized to $L(\Gamma)$.
\end{proof}

Furthermore, we have to observe that every mouse pair weaker than the generator $(K,\Psi)$ of $\Gamma$ is in $\Gamma$.
This is just the definition of a generator in the case that the former is strongly stable, but otherwise, an argument is required.

\begin{lemma}\label{905}
    Suppose that $(M,\Sigma)$ is a dropping iterate of $(K,\Psi)$.
    Then $(M,\Sigma)$ belongs to $L(\Gamma)$.
\end{lemma}
\begin{prooff}
    \item Let us assume otherwise.
    By the definition of a generator, this can only occur if $M$ is not strongly stable.

    \item Let $\Bar{M}$ be the strong $k(M)$-core of $M$ and let $\Bar{\Sigma}$ be the pullback of $\Sigma$ along the anti-core embedding.

    \item Let $D$ be the order zero measure on $\eta^M_{k(M)}$ \footnote{
    We remind the reader that $\eta^M_{k(M)}$ denotes the $r\Sigma_{k(M)}^M$-cofinality of $\rho_{k(M)}^M$.
    } and let
    $$(N,\Tau):=\Ult((\Bar{M},\Bar{\Sigma}),D).$$
    Premouse $N$ is strongly stable.

    \item Since $(N,\Tau)\not\in L(\Gamma)$, it follows that $(K,\Psi)\leq^* (N,\Tau)$.
    Thus, there exist a stack $t$ on $(N,\Tau)$ and a nearly elementary embedding 
    $$\sigma:(K,\Psi)\to(\M_\infty(t),\Tau_t).$$

    \item Let $\Bar{s}$ be the stack $(D)^\frown t$ on $(\Bar{M},\Bar{\Sigma})$, let $s$ be the lift of $\Bar{s}$ via the anti-core embedding $\Bar{M}\to M$, and let 
    $$\pi:\M_\infty(t)\to P\unlhd\M_\infty(s)$$
    be the embedding from lifting.
    We now have a nearly elementary embedding
    $$\pi\circ\sigma:(K,\Psi)\to (P,\Sigma_{s,P}).$$

    \item On the other hand, $(P,\Sigma_{s,P})$ is a dropping iterate of $(K,\Psi)$: we first have a dropping stack from $(K,\Psi)$ to $(M,\Sigma)$ and then the stack $s^\frown(P)$.
    This contradicts Dodd-Jensen Lemma\footnote{
    Cf. \cite[Theorem 9.3.4]{steel2023comparison}.
    }.
\end{prooff}

Having done these technical preparations, we are ready to make a very important point.
Namely, if $(M,\Sigma)$ is a generator for a non-limit Solovay pointclass of $L(\Gamma)$, then $(K,\Psi)$ is of course stronger than $(M,\Sigma)$, but it turns out that in a minimal comparison witnessing this, neither tree drops.

\begin{proposition}\label{min comp no drop}
    Suppose that
    \begin{assume}
        \item $\alpha<\Omega^{L(\Gamma)}$,
        
        \item $(M,\Sigma)$ is a generator for some $\Delta_{\theta_{\alpha+1}}^{L(\Gamma)}$,

        \item $(\T,\uu)$ is a minimal comparison of $(M,\Sigma)$ with $(K,\Psi)$.
    \end{assume}
    Then the main branch of $\uu$ does not drop.\footnote{
    The proof of this proposition follows closely the proof of \cite[Lemma 4.4]{steel2023mouse}, but unfortunately, we cannot directly use that lemma.
    }
\end{proposition}
\begin{prooff}
    \item Let us assume otherwise and let $(M',\Sigma')$ and $(K',\Psi')$ be the last pairs of $\T$ and $\uu$, respectively.
    We have that $(M',\Sigma')\unlhd^*(K',\Psi')$ and that $\Psi'\in L(\Gamma)$ (cf. Lemma \ref{905}).
    
    \item\label{630}\claim For all $(N,\Tau)$ for which
    $$(M',\Sigma')\unlhd (N,\Tau)\unlhd (K',\Psi'),$$
    for all $\eta<o(M')$, for all $r\in N$, we have that the set
    $$\hull^N_{k(N)+1}(\eta\cup\{r\})$$
    is not cofinal below $o(M')$.\footnote{
    We are essentially claiming that $o(M')$ is $r\Sigma_{k(N)}$-regular in $N$.
    }
    \begin{prooff}
        \item Let us assume otherwise and let $N$ be the counter-example with the least possible $l(N)$.

        \item By the minimality of $N$, we get that $o(M')$ is $r\Sigma^N_{k(N)}$-regular in $N$.
        (Note that the case $N=M'$ is covered by Lemma \ref{generator Solovay successor}.)

        \item This regularity together with the fact that $(M',\Sigma')\unlhd^*(N,\Tau)$ implies that
        \begin{parts}
            \item $\M_\infty(M',\Sigma')\unlhd^*\M_\infty(N,\Tau)$,

            \item $\pi_{(N,\Tau),\infty}\rest M'=\pi_{(M',\Sigma'),\infty}$,

            \item $\pi_{(N,\Tau),\infty}(o(M'))=o(\M_\infty(M',\Sigma'))$.
        \end{parts}

        \item\label{653} Let $\eta<o(M')$ and $r\in N$ be such that the set 
        $$\hull^N_{k(N)+1}(\eta\cup\{r\})$$
        is cofinal below $o(M')$ and let 
        $$(\eta_\infty,r_\infty):=\pi_{(N,\Tau),\infty}((\eta,r)).$$
        We have that the set
        $$\hull^{\M_\infty(N,\Tau)}_{k(N)+1}(\eta_\infty\cup\{r_\infty\})$$
        is cofinal below $o(\M_\infty(M',\Sigma'))$.\footnote{
        The fact that we have enough elementarity to push this statement about the cofinality from $N$ to the mouse limit follows from the pfs version of \cite[Corollary 4.3]{mitchell1994fine}.
        }

        \item Since $\Tau\in\Gamma$, we get that $\M_\infty(N,\Tau)\in\hh_\Gamma$.
        Thus, 
        $$\cof^{\hh_\Gamma}(o(\M_\infty(M',\Sigma')))\leq\eta_\infty<o(\M_\infty(M',\Sigma')).$$

        \item However, Lemma \ref{generator Solovay successor} shows that $o(\M_\infty(M',\Sigma'))=\theta_{\alpha+1}$, so it is Woodin in $\hh_\Gamma$.
        This is a contradiction.
    \end{prooff}

    \item\claim There exists a \underline{limit} ordinal $\gamma$ such that $\len(\uu)=\gamma+1$.
    \begin{prooff}
        \item Let us assume otherwise.
        Then there are $\xi$ and $\eta$
        such that $\len(\uu)=\eta+2$ and $\xi$ is the $\uu$-predecessor of $\eta+1$.
        

        \item Recall that 
        $$K'=\M^\uu_{\eta+1}=\Ult(\M^{*\uu}_{\eta+1},E^\uu_\eta),$$
        where $\M^{*\uu}_{\eta+1}$ is the longest initial segment to which $E^\uu_\eta$ can be applied.

        \item Since the main branch of $\uu$ drops, it follows that $K'$ is sound above $\lambda(E^\uu_\eta)$.

        \item In other words, the equality
        $$K'=\hull^{K'}_{k(K')+1}(\lambda(E^\uu_\eta)\cup\{p(K')\})$$
        holds.
        Since Claim \ref{630} excludes this possibility in the case that $\lambda(E^\uu_\eta)< o(M')$, we conclude that 
        $$\lambda(E^\uu_\eta)\geq o(M').$$

        \item Let us denote by $\kappa$ the critical point of $E^\uu_\eta$.
        The minimality of $(\T,\uu)$ implies that $\kappa\leq o(M')$.
        
        \item Thus, we know so far that
        $$\kappa\leq o(M')\leq\lambda(E^\uu_\eta)<\len(E^\uu_\eta).$$
        The next thing we verify is that the first inequality is strict.

        \item \textbf{Subclaim.} $\kappa<o(M')$
        \begin{proof}
            Let us assume otherwise, i.e. that $\kappa=o(M')$.
            Theorem \ref{generator Solovay successor} implies that 
            $$o^{M'}(\beta^{M'})=\tau^{M'}=o(M')=\kappa.$$
            Since $M'\lhd K'$ and $E^\uu_\eta$ is on the $K'$-sequence and overlaps $o(M')$, it follows that
            $$o^{K'}(\beta^{M'})>\kappa,$$
            but this contradicts the fact that $M'\unlhd^* K'$, i.e. that $o(M')$ is a cutpoint of $K'$.
        \end{proof}

        \item Thus, $\kappa<o(M')<\len(E^\uu_\eta)$.
        This in particular means that $\kappa$ is a cardinal in $M'$.
        
        \item\label{781} \textbf{Subclaim.} $\kappa$ is a strong cutpoint in $M'$.
        \begin{proof}
            Let us assume otherwise, i.e. that for some $\alpha<\kappa$, we have that $o^{M'}(\alpha)\geq\kappa$.
            This implies that 
            $$o^{K'}(\alpha)\geq\lambda(E^\uu_\eta)>o(M'),$$
            which contradicts $M'\unlhd^* K'$.
        \end{proof}

        \item Subclaim \ref{781} ensures that $\kappa\leq\beta^{M'}$.

        \item \textbf{Subclaim.} $\rho(K')\leq\kappa$
        \begin{prooff}
            \item If $\eta+1\in D^\uu$, then $\rho(K')=\rho(\M^{\uu*}_{\eta+1})\leq\kappa$ by definition.

            \item Otherwise, if $\eta+1\not\in D^\uu$, we have that
            $$[0,\xi]_\uu\cap D^\uu\not=\emptyset.$$
            Let $\alpha+1$ be the largest element of this set.

            \item It follows that
            $$\rho(K')=\rho(\M_{\alpha+1}^{\uu*})\leq\crit(E^\uu_{\alpha+1})<\lambda(E^\uu_{\alpha+1})\leq\crit(E^\uu_{\eta+1})=\kappa.$$

            \item In either case, we get that $\rho(K')\leq\kappa$.
        \end{prooff}

        \item Thus, $\rho(K')\leq\kappa\leq\beta^{M'}$, which contradicts Lemma \ref{fullness of Solovay generators}.
        (One should recall the first line of the proof where we assumed that $(K',\Psi')$ is a dropping iterate of $(K,\Psi)$, which meant that $\Psi'\in L(\Gamma)$.)
    \end{prooff}

    \item We have just argued that $\len(\uu)=\gamma+1$ for some limit $\gamma$.
    Let us now pick $\alpha$ on the main branch of $\uu$ so that there are no drops in the interval $[\alpha,\gamma)_\uu$ and so that there is an $i^\uu_{\alpha,\gamma}$-preimage $\delta$ of $o(M')$.

    \item Let $F$ be the extender applied at $\alpha$ along $[0,\gamma)_\uu$.
    We have that $\delta\geq\crit(F)$.

    \item Let $\alpha^*$ be the successor of $\alpha$ in $[0,\gamma)_\uu$, let $\delta^*:=i_{\alpha,\alpha^*}^\uu(\delta)$, let $\eta$ be the supremum of generators of $F$, and let $W:=\M_{\alpha^*}^\uu$.
    We have that $\eta<\lambda_F\leq\delta^*$ \footnote{
    We are using here the assumption that there are no mouse pairs with superstrong cardinals, which we made at the beginning of the section.
    }
    and that
    $$W=\hull^W_{k(W)+1}(\eta\cup\{p(W)\}).$$

    \item In particular, $\hull^W_{k(W)+1}(\eta\cup\{p(W)\})$ is cofinal below $\delta^*$.
    
    \item By the elementarity\footnote{
    The level of the elementarity that we need again follows from the pfs version of \cite[Corollary 4.3]{mitchell1994fine}.
    } and the continuity\footnote{
    We are using here that $o(M')$ is not measurable and is $r\Sigma_{k(K')}$-regular in $K'$.
    } at $\delta^*$ of the embedding $i_{\alpha^*,\gamma}^\uu:W\to K'$, it follows that 
    $$\hull^{K'}_{k(K')+1}(i^\uu_{\alpha^*,\gamma}(\eta)\cup\{p(K')\})$$
    is cofinal below $i_{\alpha^*,\gamma}^\uu(\delta^*)=o(M')$, which is in contradiction with Claim \ref{630}.
\end{prooff}

What we have just proven essentially says that we may assume that a generator $(M,\Sigma)$ for a non-limit Solovay pointclass of $L(\Gamma)$ is an initial segment of $(K,\Psi)$.
The next point that we want to make is that $K$ does not project strictly across the height of $M$.

\begin{proposition}\label{min comp proj}
    Suppose that
    \begin{assume}
        \item $\alpha<\Omega^{L(\Gamma)}$,
        
        \item $(M,\Sigma)$ is a generator for some $\Delta_{\theta_{\alpha+1}}^{L(\Gamma)}$,

        \item $(\T,\uu)$ is a minimal comparison of $(M,\Sigma)$ with $(K,\Psi)$.
    \end{assume}
    Then $o(\M_\infty^\T)\leq \rho^-(\M_\infty^\uu)$.\footnote{
    The proof is similar to the proof of Claim 1 of the proof of \cite[Lemma 4.4]{steel2023mouse}.
    }
\end{proposition}
\begin{prooff}
    \item Let us assume otherwise and let $(M',\Sigma')$ and $(K',\Psi')$ be the final models of $\T$ and $\uu$, respectively.
    The situation is that
    \begin{parts}
        \item neither $\T$ nor $\uu$ drops,

        \item $(M',\Sigma')\unlhd^* (K',\Psi')$,

        \item $\rho^-(K')<o(M')$.
    \end{parts}

    \item Let us denote by $(P,\Upsilon)$ the shortest initial segment of $(K',\Psi')$ for which there exist $\eta<o(M')$ and $r\in P$ so that
    $$\Hull^P_{k(P)+1}(\eta\cup\{r\})$$
    is cofinal in $o(M')$.\footnote{
    We are essentially picking the least initial segment that singularizes $o(M')$ at the \underline{next} degree of the fine structure.
    Such $P$ exists since this singularization occurs for $(K')^-$.
    }

    \item Observe that
    \begin{parts}
        \item $(M',\Sigma')\unlhd (P,\Upsilon)\lhd (K',\Psi')$,

        \item $o(M')\leq\rho^-(P)$.
    \end{parts}


    \item \textbf{Case I.} $P$ is strongly stable.
    \begin{prooff}
        \item This assumption implies that $\M_\infty(P,\Upsilon)$ belongs to $\HOD$ of $L(\Gamma)$, i.e. $\M_\infty(P,\Upsilon)\in\hh_\Gamma$.\footnote{
        Cf. \cite[Proposition 2.1]{steel2023mouse}.
        }

        \item Considering that $o(M')$ is a clean cutpoint of $P$ and is $r\Sigma^P_{k(P)}$-regular, we get that
        $$\M_\infty(M',\Sigma')\unlhd\M_\infty(P,\Upsilon)$$
        and that 
        $$\pi_{(P,\Upsilon),\infty}(o(M'))=o(\M_\infty(M',\Sigma'))=\theta_{\alpha+1}^{L(\Gamma)}.$$

        \item For all $\eta<o(M')$ and $r\in P$ for which
        $$\Hull^P_{k(P)+1}(\eta\cup\{r\})$$
        is cofinal in $o(M')$, we have that
        $$\Hull^{\M_\infty(P,\Upsilon)}_{k(P)+1}(\pi_{(P,\Upsilon),\infty}(\eta)\cup\{\pi_{(P,\Upsilon),\infty}(r)\})$$
        is cofinal in $\theta_{\alpha+1}^{L(\Gamma)}$.

        \item However, $\M_\infty(P,\Upsilon)\in\hh_\Gamma$, so it follows that $\theta_{\alpha+1}^{L(\Gamma)}$ is not regular in $\hh_\Gamma$, which of course absurd.
    \end{prooff}

    \item \textbf{Case II.} $P$ is not strongly stable.
    \begin{prooff}
        \item Let $\Bar{P}$ be the $k(P)^\mathrm{th}$ strong core of $P$ and let $\Bar{\Upsilon}$ be the pullback strategy of $\Upsilon$ along the anti-core embedding.
        The fact that $o(M')\leq\rho^-(P)$ yields that $(M',\Sigma')\unlhd (\Bar{P},\Bar{\Upsilon})$.

        \item Let $D$ be the order zero measure on $\eta^P_{k(P)}$ on the $\Bar{P}$-sequence and let
        $$(P^*,\Upsilon^*):=(\Ult(\Bar{P},D),\Bar{\Upsilon}_D).$$
        
        \item Since $o(M')$ is a clean cutpoint of $P$, it follows that either $D$ is on the $M'$-sequence or that the critical point of $D$ is strictly above $o(M')$.
        This means that in the respective cases, either\footnote{
        For the first case, one should recall that $o(M')$ is $r\Sigma^{\Bar{P}}_{k(P)}$-regular.
        } 
        $$i^{\Bar{P}}_D(M')=\Ult(M',D)\mbox{ or }i^{\Bar{P}}_D(M')=M'.$$

        \item Again in the corresponding cases, let us denote by $(M^*,\Sigma^*)$ either the pair $(\Ult(M',D), \Sigma'_D)$ or the pair $(M',\Sigma')$.
        The situation now is as follows:
        \begin{parts}
            \item $(P^*,\Upsilon^*)$ is strongly stable,
            
            \item $(M^*,\Sigma^*)\unlhd^* (P^*,\Upsilon^*)$,

            \item $(M^*,\Sigma^*)$ is a generator for $\Delta_{\theta_{\alpha+1}}^{L(\Gamma)}$,

            \item $o(M^*)$ is $r\Sigma^{P^*}_{k(P^*)}$-regular,

            \item there exist $\eta<o(M^*)$ and $r\in P^*$ such that
            $$\Hull^{P^*}_{k(P^*)+1}(\eta\cup\{r\})$$
            is cofinal in $o(M^*)$.
        \end{parts}

        \item This suffices to run the argument of Case I and obtain a contradiction in the same way.
    \end{prooff}

    \item We showed that both cases lead to a contradiction, which concludes the proof.
\end{prooff}

We can now conclude that the direct limit associated to $(K,\Psi)$ extends $\hh_\Gamma$ and does not project strictly across its height.

\begin{corollary}\label{1012}
    It holds that $\hh_\Gamma\unlhd\M_\infty(K,\Psi)$ and that 
    $$\rho^-(\M_\infty(K,\Psi))\geq\Theta_\Gamma.$$
\end{corollary}
\begin{prooff}
    \item Let $\alpha<\Omega^{L(\Gamma)}$ be arbitrary.
    We want to show that $\hh_\Gamma\para\theta_{\alpha+1}^{L(\Gamma)}\unlhd\M_\infty(K,\Psi)$ and that $\rho^-(\M_\infty(K,\Psi))>\theta_{\alpha+1}^{L(\Gamma)}$.

    \item\label{813} Let $(M,\Sigma)$ be a generator for $\Delta_{\theta_{\alpha+1}}^{L(\Gamma)}$ and let $(N,\Tau)$ be a generator for $\Delta_{\theta_{\alpha+2}}^{L(\Gamma)}$.
    Without loss of generality, we may assume that $(M,\Sigma)\lhd^* (N,\Tau)$.

    \item\label{816} By \cite[Lemma 4.4]{steel2023mouse}, we have that $o(M)\leq\rho^-(N)$ and that $o(M)$ is $r\Sigma_{k(N)}$-regular in $M$.

    \item Let $(\T,\uu)$ be a minimal comparison of $(N,\Tau)$ with $(K,\Psi)$, let $(N',\Tau')$ and $(K',\Psi')$ (resp.) be the last pairs of these trees, let $M':=i^\T(M)$, and let $\Sigma':=\Tau'_{M'}$.
    Lines \ref{813} and \ref{816} imply that $(M',\Sigma')$ is a pair on the main branch of $\T$ and that $\T$ up to that pair can be seen as a normal tree on $(M,\Sigma)$.

    \item\label{821} Propositions \ref{min comp no drop} and \ref{min comp proj} imply that $\uu$ does not drop,
    $$(M',\Sigma')\lhd^*(N',\Tau')\unlhd^*(K',\Psi'),$$
    and $\rho^-(K')\geq o(N')>o(M')$.

    \item Lines \ref{816} and \ref{821} imply that $o(M')$ is $r\Sigma_{k(K')}$-regular in $K'$.
    It follows that
    $$\M_\infty(M',\Sigma')\unlhd^*\M_\infty(K',\Psi')$$
    and $\rho^-(\M_\infty(K',\Psi'))>o(\M_\infty(M',\Sigma'))$.

    \item In other words,
    $$\hh_\Gamma\para\theta_{\alpha+1}^{L(\Gamma)}=\M_\infty(M',\Sigma')\unlhd\M_\infty(K,\Psi)$$
    and $\rho^-(\M_\infty(K,\Psi))>\theta_{\alpha+1}^{L(\Gamma)}$, as required.
\end{prooff}

We are finally ready to state the description of $K$ that we require.
As pointed out in the beginning, it is necessary for us to strengthen the anti-large-cardinal hypothesis.

\begin{proposition}\label{rho- of a generator}
    Assume $\mathsf{NMLW}$.\footnote{
    That is, we are assuming that there is no mouse pair with a measurable limit of Woodin cardinals.
    }
    Then $k(K)=0$ and that $\Theta_\Gamma$ is the largest cardinal of $\M_\infty(K,\Psi)$.
\end{proposition}
\begin{prooff}
    \item It is not possible that $\rho^-(\M_\infty(K,\Psi))=\Theta_\Gamma$.
    For otherwise, the fact that $K$ is strongly stable would imply that $\pi_{(K,\Psi),\infty}\rest \rho^-(K)$ is cofinal in $\Theta_\Gamma$, contradicting the fact that $\cof(\Theta_\Gamma)>\omega$.

    \item Thus, $\rho^-(\M_\infty(K,\Psi))>\Theta_\Gamma$.
    


    \item Let us assume towards contradiction that either $k(K)>0$ or $\Theta_\Gamma$ is not the largest cardinal of $\M_\infty(K,\Psi)$.

    \item Let $(K',\Psi')$ be a non-dropping iterate of $(K,\Psi)$ such that for some $\delta<o(K')$,
    $$\pi_{(K',\Psi'),\infty}(\delta)=\Theta_\Gamma,$$
    let $P:=K'|\delta^{+K'}$, and let $\Upsilon:=\Psi'_P$.
    Our assumptions ensure that $(P,\Upsilon)\lhd (K',\Psi')$ and that $\rho^-(K')\geq o(P)$.
    (Note that we might have $\delta^{+K'}=o(K')$.)
    
    \item\label{855}\claim $(P,\Upsilon)$ is a generator for $\Gamma$.
    \begin{prooff}
        \item If suffices to show that for all $\alpha<\Omega^{L(\Gamma)}$ and for all generators $(M,\Sigma)$ for $\Delta_{\theta_{\alpha+1}}^{L(\Gamma)}$, it holds that $(M,\Sigma)\leq^*(P,\Upsilon)$.
        
        \item Let $(\T,\uu)$ be a minimal comparison of $(M,\Sigma)$ with $(K',\Psi')$, with the last pairs $(M'',\Sigma'')$ and $(K'',\Psi'')$, and let 
        $$Q:=i^\uu(P)=K''|i^\uu(\delta^{+K'}).$$
        
        \item Since $K$ is strongly stable and since $\delta$ is \underline{not} measurable\footnote{
        This is the reason why we had to strengthen our anti-large-cardinal assumption for this proposition.
        } in $K$, it follows that $(Q,\Psi''_Q)$ is a non-dropping iterate of $(P,\Upsilon)$.
        
        \item By Corollary \ref{1012}, it follows that $(M'',\Sigma'')\unlhd (Q,\Psi''_Q)$, so we get that $(M,\Sigma)\leq^*(P,\Upsilon)$.
    \end{prooff}

    \item Claim \ref{855} contradicts the fact that $(P,\Upsilon)<^*(K,\Psi)$.
\end{prooff}

\begin{corollary}
    Assume $\mathsf{NMLW}$.
    Then some ordinal is mapped to $\Theta_\Gamma$ by the direct limit embedding $\pi_{(K,\Psi),\infty}$.
    \qed
\end{corollary}

\newpage

\section{Realizability Strategy}
\label{Realizability Strategy}

\begin{declaration}
    We are assuming that
    \begin{assume}
        \item $V=L(\ps(\R))\models\AD^+$,

        \item there does \underline{not} exist a mouse pair with a measurable limit of Woodin cardinals,

        \item $\Gamma$ is a closed pointclass such that  $L(\Gamma)\models\AD_\R$ and 
        $$\omega<\cof(\Theta_\Gamma)\leq\Theta_\Gamma<\Theta,$$

        \item $(K,\Psi)$ is a generator for $\Gamma$.
        \qed
    \end{assume}
\end{declaration}

Note that these assumptions together with Theorem \ref{HPC} imply that $\HPC$ holds.
This fact will be used tacitly from now on.

\begin{notation}
    We denote by \intro{$K_\infty$} the domain of $\M_\infty(K,\Psi)$.
    \qed
\end{notation}

The goal of this section is to show that the model $L(\Gamma,K_\infty^\omega)$ has a set of reals which is not in $\Gamma$.
This set of reals is essentially the fragment of the strategy $\Psi$ acting on the non-dropping\footnote{
Recall that the non-dropping trees are exactly those which (1) have the final model and (2) the branch leading to the final model does not drop.
} normal trees.
In order to see that the model $L(\Gamma,K_\infty^\omega)$ can reconstruct this fragment, we will have to give a sufficiently absolute description of it.
The idea how to proceed comes from \cite[Section 6.3]{sargsyan2015hod}.
Roughly, we have enough information in $\Gamma$ to iterate initial segments of $K$, while for picking a non-dropping branch through a tree, we can ask whether the model associated to a branch embeds into $K_\infty$ in a commutative manner (i.e. whether the branch is ``realizable'').
In the later case, the point is that all the embeddings belong to $K_\infty^\omega$.

\begin{notation}
    Let $M$ be a premouse which is elementarily equivalent\footnote{
    We mean that $k(M)=k(K)$ and that these premice have the same $r\Sigma_{k(K)}$-theory.
    } to $K$.
    Then the ordinals \intro{$\lambda^M$} and \intro{$\delta_\alpha^M$} (for $\alpha\leq{\lambda}^M$) are defined in such a way that the sequence $(\delta_\alpha^M : \alpha\leq {\lambda}^M)$ is a continuous, strictly increasing enumeration of cutpoint Woodin cardinals and their limits in $M$.
    We also set \intro{$\delta^M$} to be $\delta^M_{\lambda^M}$ and \intro{$\delta^M_{-1}$} to be 0.
    \qed
\end{notation}

By the results of Section \ref{Description of a Generator}, $\lambda^M$ is a limit ordinal, $\delta^M$ is the largest cardinal of $M$, and $k(M)=0$.
Furthermore, under the direct limit embedding
$$\pi_{(K,\Psi),\infty}:K\to\M_\infty(K,\Psi),$$
the ordinal $\delta^K$ is mapped to $\Theta_\Gamma$, $\M_\infty(K,\Psi)\para\Theta_\Gamma=\hh_\Gamma$, and the ordinals $\delta_\alpha^K$ (for $\alpha<\lambda^M$) are mapped to members of the Solovay sequence of $L(\Gamma)$.

The sequence $(\delta_\alpha^K : \alpha\leq \lambda^M)$ naturally divides $K$ and its iterates into windows.
While using the extenders between two successive members of this sequence, we only need to know the iteration strategy for the initial segment up to the larger member.
This is what we think of as a window.
The fragment of the iteration strategy for a single window is contained in $\Gamma$, but the complexity increases when changing a window.
In that case, we pick the correct branch of the tree by asking whether it can be realized into $K_\infty$.
This is an intuitive explanation why the information in $\Gamma\cup K_\infty^\omega$ is enough to reconstruct the fragment of $\Psi$ concerning non-dropping trees.

We will know work on formalizing the intuition just given.
The first step is to make sure that we can choose canonically the strategy for any single window.
This will be enough for the part of the previous description concerning iterations that stay inside one window.

\begin{lemma}
    Suppose that
    \begin{assume}

        \item $M$ is a countable premouse elementarily equivalent to $K$,

        \item $\sigma:\lfloor M\rfloor\xrightarrow[\Sigma_1]{\ \mathrm{cof}\ }K_\infty$,\footnote{
        We remind the reader that $\lfloor M\rfloor$ denotes the domain of the premouse $M$.
        }

        \item $\hh_\Gamma\in\ran(\sigma)$,

        \item $\xi<\lambda^M$.
    \end{assume}
    Then there exists at most one strategy $\Upsilon$ such that $(M\para\delta^M_{\xi+1},\Upsilon)$ is a generator for $\Delta_{\sigma(\delta_{\xi+1}^M)}^{L(\Gamma)}$ and $\sigma\rest(M\para\delta^M_{\xi+1})=\pi_{(M\para\delta^M_{\xi+1},\Upsilon),\infty}$.
\end{lemma}
\begin{proof}
    This is immediate from \cite[Proposition 2.1(b)]{steel2023comparison}.
\end{proof}

\begin{definition}
    Suppose that
    \begin{assume}
        \item $M$ is a countable premouse elementarily equivalent to $K$,

        \item $\sigma:\lfloor M\rfloor\xrightarrow[\Sigma_1]{\ \mathrm{cof}\ }K_\infty$,

        \item $\hh_\Gamma\in\ran(\sigma)$.
    \end{assume}
    We say that $\sigma$ is \intro{$\Gamma$-certified} iff for all $\xi<\lambda^M$, there exists a strategy $\Upsilon$ such that $(M\para\delta^M_{\xi+1},\Upsilon)$ is a generator for $\Delta_{\sigma(\delta_{\xi+1}^M)}^{L(\Gamma)}$ and $\sigma\rest(M\para\delta^M_{\xi+1})=\pi_{(M\para\delta^M_{\xi+1},\Upsilon),\infty}$.
    \qed
\end{definition}

If the premouse $M$ of the previous definition is a non-dropping iterate of $(K,\Psi)$, it will be $\Gamma$-certified.
The previous lemma and definition together mean that a $\Gamma$-certified premouse $M$ has a unique strategy with the relevant properties for each initial segment $M\para\delta^M_{\xi+1}$.

\begin{definition}
    Suppose that
    \begin{assume}

        \item $M$ is a countable premouse elementarily equivalent to $K$,

        \item $\sigma:\lfloor M\rfloor\to K_\infty$ is $\Gamma$-certified,

        \item $\delta$ is a cutpoint Woodin cardinal of $M$.
    \end{assume}
    Then we denote by \intro{$\Upsilon^\sigma_{\delta}$} the unique strategy $\Upsilon$ such that $(M\para\delta,\Upsilon)$ is a generator for $\Delta_{\sigma(\delta)}^{L(\Gamma)}$ and $\sigma\rest(M\para\delta)=\pi_{(M\para\delta,\Upsilon),\infty}$.
    \qed
\end{definition}

Observe that the model $L(\Gamma, \sigma)$ can correctly identify each $\Upsilon^\sigma_\delta$, which means that we can simply use this strategy to iterate below $\delta$.
This was the case of iterations within a single window.
We will now proceed to study the case when windows change.
If we have built a normal iteration tree with a final model $M$ and then picked an extender $E$ outside of the current window, this means that the branch leading to $M$ does not drop.
Let us distinguish nodes where a change of the window can occur.

\newcommand{\Ter}{\mathsf{Ter}}
\begin{definition}
    Suppose that $\T$ is a countable normal tree on $K$.
    Then for $\alpha<\len(\T)$, we say that $\alpha$ a \intro{terminal node} of $\T$ iff the branch $[0,\alpha]_\T$ does not drop.
    We denote by \intro{$\Ter(\T)$} the set of all terminal nodes of $\T$.
    \qed
\end{definition}

If we have built a normal tree whose last node is terminal, then the tree can easily be continued in a way that changes the window.
In particular, in order to go forward, we might require an iteration strategy for a larger initial segment.
These strategies are derived from the realization into $K_\infty$, so we want to make sure to record the realizations of the terminal nodes as we build an iteration.

\begin{definition}\label{rs-realization}
    Suppose that $\T$ is a countable normal tree on $K$.
    Then a \intro{$\pi_{(K,\Psi),\infty}$-realization} of $\T$ is a sequence
    $$\Vec{\sigma}=(\sigma_\alpha : \alpha\in\Ter(\T)),$$
    such that for all $\alpha,\beta\in\Ter(\T)$, for all $\gamma\in (\alpha,\len(\T))$, it holds that
    \begin{parts}
        \item\label{rs-212} $\sigma_\alpha:\lfloor\M^\T_\alpha\rfloor\to K_\infty$ is $\Gamma$-certified,

        \item\label{rs-110} if $\alpha=0$, then $\sigma_\alpha=\pi_{(K,\Psi),\infty}$,

        \item\label{rs-216} if $\alpha\leq_\T\beta$, then the diagram
        \begin{center}
            \begin{tikzcd}
                K_\infty & \\
                \lfloor\M^\T_\alpha\rfloor \arrow[u, "\sigma_\alpha"] \arrow[r, "i^\T_{\alpha,\beta}" swap] & \lfloor\M^\T_\beta\rfloor \arrow[ul, "\sigma_\beta" swap]
            \end{tikzcd}
        \end{center}
        commutes,

        \item\label{rs-121} if $\alpha+1<\len(\T)$, then $\len(E^\T_\alpha)<\delta^{\M^\T_\alpha}$,

        \item\label{rs-123} for all $\xi<\lambda^{\M^\T_\alpha}$, letting $\delta:=\delta_{\xi+1}(\M^\T_\alpha)$, if $\T\rest [\alpha,\gamma]$ is a tree on $\M^\T_\alpha\para\delta$, then
        \begin{parts}
            \item $\T\rest [\alpha,\gamma]$ is according to $\Upsilon^{\sigma_\alpha}_\delta$, and

            \item if $\gamma$ is a terminal node of $\T$, then, letting $\delta^*:=i^\T_{\alpha,\gamma}(\delta)$, the last pair of $\T\rest [\alpha,\gamma]$, seen as a tree on $(\M^\T_\alpha\para\delta,\Upsilon^{\sigma_\alpha}_\delta)$, is the pair $(\M^\T_\gamma\para\delta^*,\Upsilon^{\sigma_\gamma}_{\delta^*})$.
            \qed
        \end{parts}


    \end{parts}
\end{definition}

\begin{definition}
    Suppose that $\T$ is a countable normal tree on $K$.
    Then $\T$ is \intro{$\pi_{(K,\Psi),\infty}$-realizable} iff there exists a $\pi_{(K,\Psi),\infty}$-realization of $\T$.
    \qed
\end{definition}

We will be able to identify the restriction of $\Psi$ to the non-dropping trees by characterizing it as the strategy that produces $\pi_{(K,\Psi),\infty}$-realizable trees.
The crucial point is that this characterization is absolute between $V$ and $L(\Gamma,K_\infty^\omega)$.\footnote{
Note that the model $L(\Gamma,K_\infty^\omega)$ does not need to compute the embedding $\pi_{(K,\Psi),\infty}$ since it can use it as a parameter.
}

\begin{lemma}
    Suppose that
    \begin{assume}

        \item $\T$ is a countable, normal tree on $K$,

        \item $\Vec{\sigma}$ is a $\pi_{(K,\Psi),\infty}$-realization of $\T$.
    \end{assume}
    Then 
    \begin{parts}
        \item $\T$ is according to $\Psi$, and

        \item for all $\alpha\in\Ter(\T)$, we have that $\sigma_\alpha=\pi_{(\M^\T_\alpha,\Psi_{\T\rest [0,\alpha]}),\infty}$.
    \end{parts}
\end{lemma}
\begin{prooff}
    \item We show this by induction on the length of $\T$.
    To that end, let us suppose that
    \begin{parts}
        \item $\T$ is of the length $\gamma+1$, where $\gamma$ is a limit ordinal,

        \item that $\T\rest\gamma$ is according to $\Psi$,

        \item for all $\alpha\in\Ter(\T)\cap\gamma$, we have that $\sigma_\alpha=\pi_{(\M^\T_\alpha,\Psi_{\T\rest [0,\alpha]}),\infty}$.
    \end{parts}
    We want to show that $[0,\gamma)_\T=\Psi(\T\rest\gamma)$ and that if $\gamma\in\Ter(\T)$, then $\sigma_\gamma=\pi_{(\M^\T_\gamma,\Psi_{\T}),\infty}$.

    \item The set $W\subseteq\Ter(\T)\cap\gamma$ and the sequence $(\delta_\alpha : \alpha\in W)$ are defined recursively so that for all $\alpha<\gamma$:
    \begin{parts}

        \item if $\alpha\in W$, then $\delta_\alpha=\delta^{\M^\T_\alpha}_{\xi+1}$, where $\xi<\lambda^{\M^\T_\alpha}$ is the least such that $\len(E^\T_\alpha)<\delta^{\M^\T_\alpha}_{\xi+1}$,

        \item $\alpha\in W$ if and only if for all $\beta\in W\cap\alpha$, $\T\rest [\beta,\alpha+1]$ is \underline{not} a tree on $\M^\T_\beta\para\delta_\beta$.
    \end{parts}
    Observe that $0\in W$ and that $W$ is closed in $\gamma$.

    \item \textbf{Case I.} $W$ is cofinal in $\gamma$.
    \begin{prooff}
        \item In this case, there is a unique cofinal branch through $\T$: the one obtain by closing $W$ downwards for $\leq_\T$.
        This means that necessarily $[0,\gamma)_\T=\Psi(\T\rest\gamma)$ and that this branch does not drop.

        \item It remains to verify that $\sigma_\gamma=\pi_{(\M^\T_\gamma,\Psi_{\T}),\infty}$.

        \item Let $x\in\M^\T_\gamma$ be arbitrary.
        There exists $\xi\in W$ and $\Bar{x}\in\M^\T_\xi$ such that $i^\T_{\xi,\gamma}(\Bar{x})=x$.

        \item We then have that
        \begin{eqnarray}
            \sigma_\gamma(x) & = & \sigma_\gamma(i^\T_{\xi,\gamma}(\Bar{x})) \label{rs-303} \\
            & = & \sigma_\xi(\Bar{x}) \label{rs-304} \\
            & = & \pi_{(\M^\T_\xi,\Psi_{\T\rest [0,\xi]}),\infty}(\Bar{x}) \label{rs-305} \\
            & = & \pi_{(\M^\T_\gamma,\Psi_\T),\infty}(i^\T_{\xi,\gamma}(\Bar{x})) \label{rs-306} \\
            & = & \pi_{(\M^\T_\gamma,\Psi_\T),\infty}(x),
        \end{eqnarray}
        where
        \begin{itemize}
            \item[(\ref{rs-304})] follows from property \ref{rs-216} of Definition \ref{rs-realization} of a $\pi_{(K,\Psi),\infty}$-realization,

            \item[(\ref{rs-305})] follows by the inductive hypothesis,

            \item[(\ref{rs-306})] follows from the fact that $\T\rest [\xi,\gamma]$ is a normal tree on 
            $$(\M^\T_\xi,\Psi_{\T\rest [0,\xi]}).$$
        \end{itemize}
    \end{prooff}

    \item \textbf{Case II.} $W$ is bounded in $\gamma$.
    \begin{prooff}
        \item Since $W$ is closed in $\gamma$, there exists the largest element $\alpha$ of $W$.
        
        \item We have that $\T\rest [\alpha,\gamma]$ is a tree on $(\M^\T_\alpha\para\delta_\alpha,\Upsilon^{\sigma_\alpha}_{\delta_\alpha})$.

        \item\label{1286} Since $\sigma_\alpha=\pi_{(\M^\T_\alpha,\Psi_{\T\rest [0,\alpha]}),\infty}$, we have that $\Upsilon^{\sigma_\alpha}_{\delta_\alpha}=\Psi_{\T\rest [0,\alpha],\M^\T_\alpha\para\delta_\alpha}$.\footnote{
        In other words, we are saying that the pullback strategy for $\M^\T_\alpha\para\delta_\alpha$ along $\sigma_\alpha$ is exactly the same as the tail strategy of $\Psi$ obtained by first iterating along the tree $\T\rest [0,\alpha]$ and then dropping to $\M^\T_\alpha\para\delta_\alpha$.
        }

        \item\label{rs-224} By the property \ref{rs-123} of Definition \ref{rs-realization}, it holds that $\T\rest [\alpha,\gamma]$ is a tree on $(\M^\T_\alpha\para\delta_\alpha,\Upsilon^{\sigma_\alpha}_{\delta_\alpha})$.

        \item Adding to this the equality of strategies in \ref{1286}, we conclude that
        $$[0,\gamma)_\T = \Upsilon^{\sigma_\alpha}_{\delta_\alpha}(\T\rest [\alpha,\gamma)) = \Psi_{\T\rest [0,\alpha],\M^\T_\alpha\para\delta_\alpha}(\T\rest [\alpha,\gamma)).$$
        

        
        \item By the internal lift consistency of $\Psi$, we get that 
        $$[0,\gamma)_\T = \Psi_{\T\rest [0,\alpha]}(\T\rest [\alpha,\gamma))=\Psi(\T\rest\gamma).$$

        \item It remains to show that if $\gamma\in\Ter(\T)$, then $\sigma_\gamma=\pi_{(\M^\T_\gamma,\Psi_\T),\infty}$.

        \item Let $\delta^*:=i^\T_{\alpha,\gamma}(\delta_\alpha)$.
        We will use the following claim.

        \item\label{rs-343} \claim $\pi_{(\M^\T_\gamma,\Psi_\T),\infty}\rest(\M^\T_\gamma\para\delta^*)=\pi_{(\M^\T_\gamma\para\delta^*,\Upsilon^{\sigma_\gamma}_{\delta^*}),\infty}$
        \begin{proof}
            Property \ref{rs-123} of Definition \ref{rs-realization} implies that $(\M^\T_\gamma\para\delta^*,\Upsilon^{\sigma_\gamma}_{\delta^*})$ is an iterate of $(\M^\T_\alpha\para\delta_\alpha,\Upsilon^{\sigma_\alpha}_{\delta_\alpha})$ via $\T\rest [\alpha,\gamma]$.
            Adding to this \ref{rs-224}, we obtain that $\Upsilon^{\sigma_\gamma}_{\delta^*}=\Psi_{\T,\M^\T_\gamma\para\delta^*}$.
            Thus,
            $$\pi_{(\M^\T_\gamma,\Psi_\T),\infty}\rest(\M^\T_\gamma\para\delta^*)=\pi_{(\M^\T_\gamma,\Psi_\T)\para\delta^*,\infty}=\pi_{(\M^\T_\gamma\para\delta^*,\Upsilon^{\sigma_\gamma}_{\delta^*}),\infty}.$$
        \end{proof}

        \item Let $x\in\M^\T_\gamma$ be arbitrary.
        Since $\T\rest [\alpha,\gamma]$ can be seen as a tree on $(\M^\T_\alpha\para\delta_\alpha,\Upsilon^{\sigma_\alpha}_{\delta_\alpha})$, we have that there exists a function $f\in\M^\T_\alpha$ and some $z\in\M^\T_\gamma\para\delta^*$ such that
        $$x=i^\T_{\alpha,\gamma}(f)(z).$$

        \item We now have that
        \begin{eqnarray}
            \sigma_\gamma(x) & = & \sigma_\gamma(i^\T_{\alpha,\gamma}(f)(z)) \\
            & = & \sigma_\gamma(i^\T_{\alpha,\gamma}(f))(\sigma_\gamma(z)) \\
            & = & \sigma_\alpha(f)(\sigma_\gamma(z)) \label{rs-348}\\
            & = & \pi_{(\M^\T_\alpha,\Psi\rest [0,\alpha]),\infty}(f)(\sigma_\gamma(z)) \label{rs-349}\\
            & = & \pi_{(\M^\T_\gamma,\Psi_\T),\infty}(i^\T_{\alpha,\gamma}(f))(\sigma_\gamma(z)) \label{rs-350}\\
            & = & \pi_{(\M^\T_\gamma,\Psi_\T),\infty}(i^\T_{\alpha,\gamma}(f))(\pi_{(\M^\T_\gamma\para\delta^*,\Upsilon^{\sigma_\gamma}_{\delta^*}),\infty}(z)) \label{rs-351}\\
            & = & \pi_{(\M^\T_\gamma,\Psi_\T),\infty}(i^\T_{\alpha,\gamma}(f))(\pi_{(\M^\T_\gamma,\Psi_\T),\infty}(z)) \label{rs-352}\\
            & = & \pi_{(\M^\T_\gamma,\Psi_\T),\infty}(i^\T_{\alpha,\gamma}(f)(z)) \\
            & = & \pi_{(\M^\T_\gamma,\Psi_\T),\infty}(x),
        \end{eqnarray}
        where
        \begin{itemize}
            \item[(\ref{rs-348})] follows from property \ref{rs-216} of Definition \ref{rs-realization} of a $\pi_{(K,\Psi),\infty}$-realization,

            \item[(\ref{rs-349})] follows by the inductive hypothesis,

            \item[(\ref{rs-350})] follows from the fact that $\T\rest [\alpha,\gamma]$ can be seen as a tree on $(\M^\T_\alpha,\Psi_{\T\rest [0,\alpha]})$,

            \item[(\ref{rs-351})] follows from the property \ref{rs-212} of Definition \ref{rs-realization} of a $\pi_{(K,\Psi),\infty}$-realization,

            \item[(\ref{rs-352})] follows from Claim \ref{rs-343}.
        \end{itemize}

        \item This concludes the verification of Case II.
    \end{prooff}

    \item The verification of the two cases concludes the proof of the lemma.
\end{prooff}

\begin{corollary}
    Let $\T$ be a countable, non-dropping, normal tree on $K$.
    Then $\T$ is according to $\Psi$ if and only if it is $\pi_{(K,\Psi),\infty}$-realizable.
\end{corollary}
\begin{proof}
    The implication ($\Leftarrow$) follows from the lemma.
    For the other direction, let $\T$ be according to $\Psi$.
    We claim that
    $$\Vec{\sigma}:=(\pi_{(\M^\T_\alpha,\Psi_{\T\rest [0,\alpha]}),\infty} : \alpha\in\Ter(\T))$$
    is a $\pi_{(K,\Psi),\infty}$-realization of $\T$.
    Our assumption on the extenders of $K$ ensures that for all $\alpha\in\Ter(\T)$ such that $\alpha+1<\len(\T)$, we have that $\len(E^\T_\alpha)<\delta^{\M^\T_\alpha}$.
    This shows the property \ref{rs-121} of Definition \ref{rs-realization}, while properties \ref{rs-110} and \ref{rs-216} are obvious.
    The rest will follow once we verify the following.
    For all $\alpha\in\Ter(\T)$ and for all $\xi<\lambda^{\M^\T_\alpha}$, letting $\delta:=\delta_{\xi+1}(\M^\T_\alpha)$, the pair
    $$(\M^\T_\alpha\para\delta,\Psi_{\T\rest [0,\alpha],\M^\T_\alpha\para\delta})$$
    is a generator for $\Delta_{\theta_{\alpha+1}}^{L(\Gamma)}$ and
    $$\sigma\rest (\M^\T_\alpha\para\delta)=\pi_{(\M^\T_\alpha\para\delta,\Psi_{\T\rest [0,\alpha],\M^\T_\alpha\para\delta}),\infty}.$$
    The verification of this is given below.
    \begin{pfenum}
        \item Let us assume otherwise.
        Then, letting $(M,\Sigma)$ be a generator for $\Delta_{\theta_{\alpha+1}}^{L(\Gamma)}$, we have that either
        $$(\M^\T_\alpha\para\delta,\Psi_{\T\rest [0,\alpha],\M^\T_\alpha\para\delta})<^*(M,\Sigma)\mbox{ or }(M,\Sigma)<^*(\M^\T_\alpha\para\delta,\Psi_{\T\rest [0,\alpha],\M^\T_\alpha\para\delta}).$$

        \item\label{rs-303'} Let $\theta:=\pi_{(\M^\T_\alpha,\Psi),\infty}(\delta)$.
        Then $\theta$ is a successor member of the Solovay sequence of $L(\Gamma)$ and
        $$\M_\infty(\M^\T_\alpha\para\delta,\Psi_{\T\rest [0,\alpha],\M^\T_\alpha\para\delta})=\M_\infty(K,\Psi)\para\theta=\hh_\Gamma\para\theta=\M_\infty(M,\Sigma).$$

        \item If it was the case that $(\M^\T_\alpha\para\delta,\Psi_{\T\rest [0,\alpha],\M^\T_\alpha\para\delta})<^*(M,\Sigma)$, we would have that 
        $$(\M^\T_\alpha\para\delta,\Psi_{\T\rest [0,\alpha],\M^\T_\alpha\para\delta})\in\Delta_\theta^{L(\Gamma)}$$
        and consequently,
        $$o(\M_\infty(\M^\T_\alpha\para\delta,\Psi_{\T\rest [0,\alpha],\M^\T_\alpha\para\delta}))<\theta.$$
        This is clearly a contradiction.

        \item Thus, the only remaining possibility is that 
        $$(M,\Sigma)<^*(\M^\T_\alpha\para\delta,\Psi_{\T\rest [0,\alpha],\M^\T_\alpha\para\delta}).$$

        \item Let $(N,\Tau)\lhd^*(P,\Upsilon)$ be the result of a minimal comparison witnessing the above inequality.
        It follows from \cite[Lemma 4.4]{steel2023mouse} that neither of the trees in this comparison drop, that $o(N)\leq\rho^-(P)$, and that $o(N)$ is $r\Sigma^P_{k(P)}$-regular.

        \item This suffices for the conclusion that $\M_\infty(N,\Tau)\lhd\M_\infty(P,\Upsilon)$, which again contradicts \ref{rs-303'}.
    \end{pfenum}
\end{proof}

The fact just established, together with the absoluteness of the $\pi_{(K,\Psi),\infty}$-realizability, gives us the absoluteness of the relevant fragment of $\Psi$.

\begin{corollary}\label{rs-322}
    The restriction of $\Psi$ to non-dropping normal trees belongs to $L(\Gamma,K_\infty^\omega)$.
    \qed
\end{corollary}

This fragment of $\Psi$ is coded by a set of reals which cannot belong to $\Gamma$: the reason for this is the ability to define from it a surjection $\R\twoheadrightarrow\hh_\Gamma$.
The surjection in question is just the one coming from the direct limit construction for $\M_\infty(K,\Psi)$ and it witnesses that $\Theta_\Gamma$ is not the ``real'' $\Theta$.

\begin{corollary}\label{361}
    There exists a set of reals in $L(\Gamma,K_\infty^\omega)$ which is not in $\Gamma$.
\end{corollary}
\begin{prooff}

    \item Let $\ff$ consist of all non-dropping normal trees on $(K,\Psi)$.
    The mapping
    $$\ff\to\ff(K,\Psi):\T\mapsto (\M^\T_\infty,\Psi_\T)$$
    is a bijection.\footnote{
    We remind the reader that $\ff(K,\Psi)$ denotes the set of all of non-dropping iterates of $(K,\Psi)$.
    }

    \item Corollary \ref{rs-322} implies that $\ff\in L(\Gamma,K_\infty^\omega)$.

    \item For $\T,\T'\in\ff$, we define $\T\leq_\ff\T'$ iff there exists a non-dropping $\pi_{(\M^{\T},\Psi_\T),\infty}$-realizable tree $\uu$ on $\M^\T_\infty$ such that the stack $\T^\frown\uu$ normalizes to $\T'$.
    
    \item Note that for all $\T\in\ff$, $\M_\infty(\M^{\T},\Psi_\T)=K_\infty$ and consequently, 
    $$L(\Gamma,K_\infty^\omega)$$
    can correctly identify $\pi_{(\M^{\T},\Psi_\T),\infty}$-realizable trees.
    In particular, $\leq_\ff\in L(\Gamma,K_\infty^\omega)$.

    \item Observe that for all $\T,\T'\in\ff$, we have that $\T\leq_{\ff}\T'$ if and only if $(\M^{\T'}_\infty,\Psi_{\T'})$ is a non-dropping iterate of $(\M^\T_\infty,\Psi_\T)$.
    Furthermore, if $\uu$ witnesses that $\T\leq_\ff\T'$, then $\uu$ is equal to the unique normal tree on $(\M^\T_\infty,\Psi_\T)$ with the last pair $(\M^{\T'}_\infty,\Psi_{\T'})$.

    \item This means that for $\T\leq_\ff\T'$, we can define
    $$\pi_{\T,\T'}:=\pi^\uu:\M^\T_\infty\to\M^{\T'}_\infty,$$
    where $\uu$ is the unique tree witnessing that $\T\leq_\ff\T'$.

    \item In this way we obtain a family of embeddings $(\pi_{\T,\T'} : \T\leq_\ff\T')$ which belongs to $L(\Gamma,K_\infty^\omega)$.

    \item On the other hand, for all $\T\leq_\ff\T'$,
    $$\pi_{\T,\T'}=\pi_{(\M^\T_\infty,\Psi_\T),(\M^{\T'}_\infty,\Psi_{\T'})}.$$

    \item Thus, the family
    $$(\M^\T_\infty,\, \pi_{\T,\T'} : \T,\T'\in\ff,\, \T\leq_\ff\T')$$
    belongs to $L(\Gamma,K_\infty^\omega)$ and forms a directed system which is equivalent to the standard directed system leading to $\M_\infty(K,\Psi)$.

    \item Now, for each $x\in\R$ which codes a pair $(\T,\alpha)$, where $\T\in\ff$ and $\alpha<\delta^{\M^\T_\infty}$, we can map $x$ to 
    $$\pi_{\T,\infty}(x)\in K_\infty.$$
    This defines a surjection onto $\Theta_\Gamma$ inside $L(\Gamma,K_\infty^\omega)$, which concludes the proof.
\end{prooff}

\newpage

\section{Adding Sequences and Measures}
\label{Adding Sequences and Measures}

\begin{declaration}
    We are assuming that
    \begin{assume}
        \item $V=L(\ps(\R))\models\AD^+ + \AD_\R$,

        \item there does \underline{not} exist a mouse pair with a measurable limit of Woodin cardinals,

        \item $\Gamma$ is a closed pointclass such that  $L(\Gamma)\models\AD_\R$ and 
        $$\omega<\cof(\Theta_\Gamma)\leq\Theta_\Gamma<\Theta.$$
        \qed
    \end{assume}
\end{declaration}

The goal of this section is to show that for some $\kappa<\Theta$, the model 
$$L(\Gamma,\kappa^\omega)[\mu_\kappa]$$
contains a set of reals which is not in $\Gamma$.\footnote{
We remind the reader that $\mu_\kappa$ is the quasi-club filter on $[\kappa]^\omega$ introduced in Definition \ref{199}.
This filter is an ultrafilter by Proposition \ref{485}.
}
By the results of the previous section, it will suffice to find $\kappa<\Theta$ so that $K_\infty^\omega$ belongs to $L(\Gamma,\kappa^\omega)[\mu_\kappa]$, where $K_\infty$ comes from a generator $(K, \Psi)$ for $\Gamma$.
It will turn out that we can take $\kappa$ to be the ordinal height of $K_\infty$.

The first step is to show that $K_\infty$ belongs to $L(\Gamma,\kappa^\omega)[\mu_\kappa]$.
The argument that we have in mind is to give a sufficiently absolute description of this set.
We will represent $K_\infty$ as a $\mu_\kappa$-ultraproduct of the transitive collapses of the hulls of $X$ inside $\M_\infty(K,\Psi)$, where $X$ varies over the countable subsets of $\kappa$.
The point here is that this family of collapsed hulls is equal almost everywhere to a family inside $L(\Gamma,\kappa^\omega)[\mu_\kappa]$.
More concretely, almost every collapsed hull of $\M_\infty(K,\Psi)$ is a non-dropping iterate of $(K,\Psi)$ and its part below its largest cardinal is just the corresponding hull of $\hh_\Gamma$.
This means that the only remaining issue is to describe the parts of the iterates of $(K,\Psi)$ above their largest cardinals in a way that relativizes to $L(\Gamma,\kappa^\omega)[\mu_\kappa]$.

Thus, the first problem that we will address is, given a generator for $\Gamma$, compute the domain of its premouse component from its part up to its largest cardinal inside $L(\Gamma)$.
Denoting by $(K,\Psi)$ the generator, we would like to say that $K$ can be obtained by stacking the premouse components of all mouse pairs that belong to $L(\Gamma)$, extend $(K, \Psi)|\delta^K$, and project to $\delta^K$.
This corresponds to the notion of the \textit{lower part operator}, but it seems that this scenario does not fully work out in the least branch hierarchy.
For now, let us distinguish the pairs that we want to stack.

\begin{definition}
    Let $(K,\Psi)$ be a generator for $\Gamma$ and let $(N,\Tau)$ be a mouse pair in $L(\Gamma)$.
    We say that $(N,\Tau)$ \intro{anticipates} $(K,\Psi)$ iff $(K,\Psi)|\delta^K\unlhd^*(N,\Tau)$ and $\rho^-(N)=\rho_\omega(N)=\delta^K$.
    \qed
\end{definition}

We would like to say that if $(N,\Tau)$ and $(P, \Upsilon)$ both anticipate $(K, \Psi)$, then one of the premice $N$ and $P$ is an initial segment of the other and the both are initial segments of $K$, but this does not seem to be true.
Nevertheless, we will be able to show that either $N$ and $P$ have the same domain or one of them is an \textit{element} of the other.
We will also establish that both of $N$ and $P$ must be elements of $K$.
Thus, the premouse components of the anticipating pairs all align with each other and with $K$ as well, but with respect the membership relation only.
Furthermore, we will be able to argue that the domain of $K$ is the union of the definable powersets all these premice, therefor obtaining a lower-part-like characterization of the \textit{set underlining} $K$.
This will suffice for the purposes of this paper.
What follows below, excepting the final theorem of the section, is verifying what we have just elaborated.
We start from a useful technical observation.

\begin{lemma}\label{nct-26}
    Suppose that
    \begin{assume}
        \item $(K,\Psi)$ is a generator for $\Gamma$,
        
        \item $(N,\Tau)$ anticipates $(K,\Psi)$.
    \end{assume}
    Then $\Bar{\mathfrak{C}}_{k(N)}(N)=N$.\footnote{
    In other words, $N$ is strongly sound.
    }
\end{lemma}
\begin{proof}
    Otherwise, $\rho^-(N)$ would be a measurable limit of Woodin cardinals in $\Bar{\mathfrak{C}}_{k(N)}(N)$.
    This is not possible since we are assuming that there is no mouse pair with a measurable limit of Woodin cardinals.
\end{proof}

If a pair $(N, \Tau)$ anticipates $(K, \Psi)$, then $N$ and $K$ agree up to $\delta^K$, but might be incomparable above it.
It turns out that we can produce a non-dropping tree on their common part below $\delta^K$ which aligns them.
More precisely, when this tree is lifted to $N$ and $K$, the final model corresponding to $N$ is an initial segment of the final model corresponding to $K$.

\begin{lemma}\label{nct-51}
    Suppose that
    \begin{assume}
        \item $(K,\Psi)$ is a generator for $\Gamma$,
        
        \item $(N,\Tau)$ anticipates $(K,\Psi)$.
    \end{assume}
    Then there exists a non-dropping, normal tree $\T$ on $(K,\Psi)|\delta^K$ such that, when regarded as a tree on $(N,\Tau)$ and $(K,\Psi)$, it has the last models $N^*$ and $K^*$ (resp.) satisfying that $N^*\lhd K^*$.\footnote{Since $N$ projects to $\delta^K$, the lifting of a tree on $K|\delta^K$ simply means applying exactly the same extenders.
    The similar fact is true for $K$ due to $\delta^K$ being its largest cardinal.}
\end{lemma}
\begin{prooff}
    \item If $N$ is not strongly stable, then, by Lemma \ref{nct-26}, letting $D$ be the order zero measure on $\eta^N_{k(N)}$, we have that 
    $$(\Ult(N,D),\Tau_{(D)})$$
    is a strongly stable mouse pair.

    \item These facts, together with the comparison for strongly stable mouse pairs and the full normalization\footnote{
    Our reference for the full normalization is \cite{siskind2022full}.
    }, imply that there are a normal tree $\uu$ on $(N,\Tau)$ with the last pair $(N^*,\Tau^*)$ and a normal tree $\vv$ on $(K,\Psi)$ with the last pair $(K',\Psi')$ such that the main branch of $\uu$ does not drop and $(N^*,\Tau^*)\lhd^* (K',\Psi')$.
    
    \item There exists the largest ordinal $\gamma_\uu$ such that $\uu\rest\gamma_\uu$ can be regarded as a tree on $K|\delta^K$.
    There exists the analogous ordinal $\gamma_\vv$ with respect to $\vv$.

    \item We have that $\gamma_\uu$ and $\gamma_\vv$ are on the main branches of $\uu$ and $\vv$ (resp.).

    \item Let $\T_\uu$ and $\T_\vv$ be the trees $\uu\rest [0,\gamma_\uu]$ and $\vv\rest [0,\gamma_\vv]$ (resp.) regarded as trees on $(K,\Psi)|\delta^K$.
    These are non-dropping, normal trees on $(K,\Psi)|\delta^K$.

    \item All extenders $E^\uu_\xi$ for $\xi\geq\gamma_\uu$ have the critical point strictly larger than $o(\M^{\T_\uu}_{\gamma_\uu})$ and all extenders $E^\vv_\xi$ for $\xi\geq\gamma_\vv$ have the critical point strictly larger than $o(\M^{\T_\vv}_{\gamma_\vv})$.
    In both cases, the said height is a limit of Woodin cardinals and consequently, a clean cutpoint in $\M^\uu_{\gamma_\uu}$ (resp. $\M^\vv_{\gamma_\vv}$).

    \item Now, since $(N^*,\Tau^*)\lhd (K',\Psi')$, we have that one of the final pairs of $\T_\uu$ and $\T_\vv$ is an initial segment of the other.
    Since these two mouse pairs are equivalent, they must be equal.

    \item This has as a consequence that $\T_\uu=\T_\vv$.
    Letting $\T$ denote this tree and $(M^*,\Sigma^*)$ be its last pair, we have that $\T\unlhd\uu,\vv$ and
    $$(M^*,\Sigma^*)\unlhd (N^*,\Tau^*),(K',\Psi').$$

    \item Since $\uu$ is non-dropping, we must have that in fact $\uu=\T\unlhd\vv$.

    \item Let us denote by $K^*$ the model $\M^\vv_{\gamma_\vv}$, i.e. the last model of $\T$ when considered on $(K,\Psi)|\delta^K$.
    The fact that all extenders $E^\vv_\xi$ for $\xi\geq\gamma_\vv$ have the critical point strictly larger than $o(M^*)$ implies that $\ps(o(M^*))\cap K'\subseteq K^*$.

    \item Since $\rho^-(N^*)=o(M^*)$, we can conclude that $N^*\lhd K^*$.
    This shows that the tree $\T$ is as required.
\end{prooff}

Even though we do not seem to be able to pull this alignment along $\T$ back to the starting premice $N$ and $K$, we are able to pull back the alignment with respect to the membership.
We will prove this after stating another useful technical observation.

\begin{lemma}\label{one tree two lifts}
    Suppose that
    \begin{assume}
        \item $(K,\Psi)$ is a generator for $\Gamma$,
        
        \item $(N,\Tau)$ anticipates $(K,\Psi)$,

        \item $\T$ is a non-dropping normal tree on $K \para \delta^K$,

        \item $\uu$ and $\vv$ are the tree $\T$ seen as a tree on $(N,\Tau)$ and $(K,\Psi)$, respectively.
    \end{assume}
    Then for all $\alpha \leq \beta$ on the main branch of $\T$,
    \begin{parts}
        \item $\M^\T_\alpha$ is an initial segment of both $\M^\uu_\alpha$ and $\M^\vv_\alpha$,
        
        \item $\rho^-(\M^\uu_\alpha)=o(\M^\T_\alpha)$,

        \item $i^\T_{\alpha,\beta} = i^\uu_{\alpha,\beta} \restriction \M^\T_\alpha = i^\vv_{\alpha,\beta} \restriction \M^\T_\alpha$.
        \qed
    \end{parts}
\end{lemma}

\begin{proposition}\label{nct-109}
    Suppose that
    \begin{assume}
        \item $(K,\Psi)$ is a generator for $\Gamma$,
        
        \item $(N,\Tau)$ anticipates $(K,\Psi)$.
    \end{assume}
    Then $N\in K$
\end{proposition}
\begin{prooff}
    \item By Lemma \ref{nct-51}, there exists a non-dropping, normal tree $\T$ on $(K,\Psi)|\delta^K$ such that, when lifted to a tree $\uu$ on $(N,\Tau)$ with final model $\Tilde{N}$ and a tree $\vv$ on $(K,\Psi)$ with the final model $\Tilde{K}$, we have $\Tilde{N}\lhd\Tilde{K}$.
    In particular, $\Tilde{N}\in\Tilde{K}$.

    \item Let $\zeta$ on the main branch of $\T$ be the least such that
    $$\M^\uu_\zeta\in\M^\vv_\zeta.$$
    We will be done if we show that $\zeta=0$.
    
    \item Let $N^*:=\M^\uu_\zeta$ and let $K^*:=\M^\vv_\zeta$.

    \item \claim $\zeta$ is not a successor.
    \begin{prooff}
        \item Let us assume otherwise, i.e. that there exists $\gamma$ such that $\zeta=\gamma+1$.
    
        \item Let $\gamma'$ be the $\T$-predecessor of $\gamma+1$, let $N':=\M^\uu_{\gamma'}$, let $K':=\M^\vv_{\gamma'}$, and let $E:=E^\T_\gamma$.
        We have that $N'\not\in K'$ and that
        $$\Ult(N',E)=N^*\in K^*=\Ult(K',E).$$
    
        \item We will reach a contradiction by showing that $N' \in K'$.
        To that end, it suffices to show that $\mathsf{Th}^{N'}_{k(N)}(\rho^-(N')\cup\{p_{k(N)}(N')\})\in K'$.

        \item Let $\phi$ be a formula, let $s\in [\rho^-(N')]^{<\omega}$, and let
        $$N'\models\phi [s,p_{k(N)}(N')].$$
        
        \item\label{1627} It follows by the elementarity that
        $$N^*\models\phi [i^{N'}_E(s),p_{k(N)}(N^*)].$$
    
        \item There exist an $r\Sigma_{k(K)}^{K'}$-function $f:\kappa\to K'$ and $\alpha<\lambda_E$ such that $N^*=i^{K'}_E(f)(\alpha)$.
    
        \item Lemma \ref{one tree two lifts} implies that $i^{N'}_E(s)=i^{K'}_E(s)$.
        Adding to this \ref{1627}, we conclude that
        $$i^{K'}_E(f)(\alpha)\models\phi [i^{K'}_E(s),p_{k(N)}(i^{K'}_E(f)(\alpha))].$$
    
        \item It follows that for $E_\alpha$-almost every $\xi<\kappa_E$,
        $$f(\xi)\models\phi [s,p_{k(N)}(f(\xi))].$$
    
        \item Thus,
        $$\mathsf{Th}^{N'}_{k(N)}(\rho^-(N')\cup\{p_{k(N)}(N')\})$$
        consists of all those $\phi$ such that for $E_\alpha$-almost every $\xi<\kappa_E$,
        $$\phi\in\mathsf{Th}^{f(\xi)}_{k(N)}(\rho^-(N')\cup\{p_{k(N)}(f(\xi))\}).$$
    
        \item Since $E$ is close to $K'$, we get that
        $$\mathsf{Th}^{N'}_{k(N)}(\rho^-(N')\cup\{p_{k(N)}(N')\})\in K'$$
        and consequently, that $N'\in K'$, which is a contradiction.
    \end{prooff}
    
    \item \claim $\zeta$ is not a limit.
    \begin{prooff}
        \item Let us assume otherwise.
        
        \item\label{1653} Since $N^*\in K^*$, there exist $\gamma\in [0,\zeta)_\T$ and $N''\in\M^\vv_\gamma$ such that
        $$i^\vv_{\gamma,\zeta}(N'')=N^*.$$

        \item Let $N':=\M^\uu_\gamma$ and $K':=\M^\vv_\gamma$.
        We will reach a contradiction by showing that $N'\in K'$.
        
        \item\label{1659} Let $\rho:=\rho^-(N')=o(\M^\T_\gamma)$.
        Since $N''\in K'$, it suffices to show that
        $$N'=\mathsf{cHull}^{N''}_{k(N)}(\rho\cup\{p_{k(N)}(N'')\}),$$
        which in turn comes down to showing that
        $$\mathsf{Th}^{N'}_{k(N)}(\rho\cup\{p_{k(N)}(N')\}) = \mathsf{Th}^{N''}_{k(N)}(\rho\cup\{p_{k(N)}(N'')\}).$$

        \item Suppose that $\phi$ is a formula, $s\in [\rho]^{<\omega}$, and
        $$N'\models\phi [s,p_{k(N)}(N')].$$
        
        \item By the elementarity of $i^\uu_{\gamma,\zeta}$, it follows that
        $$N^*\models\phi [i^\uu_{\gamma,\zeta}(s),p_{k(N)}(N^*)].$$

        \item Lemma \ref{one tree two lifts} implies that $i^\uu_{\gamma,\zeta}(s)=i^\vv_{\gamma,\zeta}(s)$.
        Adding to this line \ref{1653}, we conclude that
        $$i^\vv_{\gamma,\zeta}(N'')\models\phi [i^\vv_{\gamma,\zeta}(s),p_{k(N)}(i^\vv_{\gamma,\zeta}(N''))].$$

        \item It follows by the elementarity of $i^\vv_{\gamma,\zeta}$ that
        $$N''\models\phi [s,p_{k(N)}(N'')].$$

        \item Since $\phi$ and $s$ were arbitrary, we conclude that
        $$\mathsf{Th}^{N'}_{k(N)}(\rho\cup\{p_{k(N)}(N')\}) = \mathsf{Th}^{N''}_{k(N)}(\rho\cup\{p_{k(N)}(N'')\}),$$
        which in turn yields $N'\in K'$, as explained in line \ref{1659}.
        We have reached a contradiction.
    \end{prooff}

    \item Thus, the only remaining possibility is that $\zeta=0$, which means that
    $$N=N^*\in K^*=K.$$
\end{prooff}

This shows that the anticipating pairs align with the generator.
Essentially the same arguments show that they also align among themselves.

\begin{proposition}\label{nct-ordering by in}
    Suppose that
    \begin{assume}
        \item $(K,\Psi)$ is a generator for $\Gamma$,

        \item $(N,\Tau)$ and $(P,\Upsilon)$ anticipate $(K,\Psi)$.
    \end{assume}
    Then either $\hat{N}\in\hat{P}$, or $\hat{N}=\hat{P}$, or $\hat{P}\in\hat{N}$.\footnote{
    Recall that $\hat{N}$ denotes the bare premouse underlining $N$, i.e. the structure obtained from $N$ by disregarding the soundness degree $k(N)$.}
\end{proposition}
\begin{proof}
    Essentially the same proof as that of Lemma \ref{nct-51} shows that there exists a non-dropping normal tree on $(K,\Psi)|\delta^K$ such that, when lifted to the tree on $(N,\Tau)$ with the final model $N^*$ and to the tree on $(P,\Upsilon)$ with the final model $P^*$, it holds that either $N^*\unlhd^*P^*$ or $P^*\unlhd^*N^*$.
    If $o(N^*)=o(P^*)$, then it follows that $\hat{N}=\hat{P}$.
    Hence, we may assume without loss of generality that $o(N^*)<o(P^*)$, i.e. that $N^*\in P^*$.
    We can now run the proof of Proposition \ref{nct-109} to establish that in fact $N\in P$.
\end{proof}

As a corollary, we obtain a lower-part-like characterization of the generator.

\begin{corollary}\label{nct-352}
    Let $(K,\Psi)$ be a generator for $\Gamma$.
    Then the domain of $K$ is equal to the union of all $J_1(N)$, where $(N,\Tau)$ anticipates $(K,\Psi)$.
\end{corollary}
\begin{prooff}
    \item If $(N,\Tau)$ anticipates $(K,\Psi)$, then it holds by Proposition \ref{nct-109} that $J_1(N)\subseteq K$.
    This verifies one inclusion.
    
    \item For the other inclusion, it suffices to show that for all $A\in\ps(\delta^K)\cap K$, there exists a mouse pair $(N,\Tau)$ which anticipates $(K,\Psi)$ so that $A\in J_1(N)$.\footnote{
    The reason why this is sufficient is that $\delta^K$ is the largest cardinal of $K$.
    }
    
    \item Since $A\in\ps(\delta^K)\cap K$, there exists $\alpha<\hat o(K)$ such that $\rho_\omega(K|\alpha)=\delta^K$ and $A\in J_1(K|\alpha)$.
    
    \item Let $k<\omega$ be such that $\rho_k(K|\alpha)=\delta^K$ and let $(N,\Tau):=(K,\Psi)|\alpha$.
    It follows from Lemma \ref{905} that $(N,\Tau)\in L(\Gamma)$, so $(N,\Tau)$ anticipates $(K,\Psi)$.
    
    \item Since $A\in J_1(N)$, this verifies the other inclusion.
\end{prooff}

We are finally ready to state the main result of the section.

\begin{theorem}\label{1882}
   Let $(K,\Psi)$ be a generator for $\Gamma$ and let $\kappa$ be the height of $\M_\infty(K,\Psi)$.
   Then $L(\Gamma,\kappa^\omega)[\mu_\kappa]$ contains a set of reals which is not in $\Gamma$.
\end{theorem}
\begin{prooff}
    \item\label{nct-374} Let us assume otherwise and let $\Mod:=L(\Gamma,\kappa^\omega)[\mu_\kappa]$.
    We have that
    \begin{parts}
        \item $\ps(\R)\cap\Mod=\Gamma$,

        \item $\kappa^\omega\subseteq\Mod$,

        \item $\mu_\kappa^\Mod=\mu_\kappa\cap\Mod$ is a supercompactness measure in $\Mod$.
    \end{parts}

    \item For $X\in [\kappa]^\omega$, let 
    $$\sigma_X:\kk_X\to\M_\infty(K,\Psi)$$
    be the anti-collapse of $\Hull^{\M_\infty(K,\Psi)}_1(X)$, where we are assuming the pullback premouse structure on $\kk_X$.

    \item \claim $\prod_{\mu_\kappa}\kk_X=\M_\infty(K,\Psi)$
    \begin{prooff}
        \item For all $X\in [\kappa]^\omega$ and for all $x\in \kk_X$, there exists the least pair $(t_{X,x},p_{X,x})$, with $t_{X,x}$ a Skolem term and $p_{X,x}\in [X]^{<\omega}$, such that
        $$\sigma_X(x)=t_{X,x}(p_{X,x}).$$

        \item For all $f\in\prod_{X\in [\kappa]^\omega}K_X$, there exists a unique pair $(t_f,p_f)$ such that for $\mu_\kappa$-almost every $X$, 
        $$(t_{X,f(X)},p_{X,f(X)})=(t_f,p_f).$$
        (This follows from the fact that $\mu_\kappa$ is a supercompactness measure.)

        \item Let $\kk$ be the ultraproduct $\prod_{\mu_\kappa\cap\Mod}\kk_X$, as computed inside $\Mod$.
        The mapping
        $$j:\kk\longrightarrow \M_\infty(K,\Psi) : [f]\mapsto t_f(p_f)$$
        is well defined and elementary.

        \item We want to see that $j$ is a surjection.
        Let $y\in\M_\infty(K,\Psi)$ be arbitrary and let us find $f\in\prod_{X\in [\kappa]^\omega}K_X\cap\Mod$ such that $j([f])=y$.

        \item For all $X\in [\kappa]^\omega$ satisfying that $y\in\Hull^{\M_\infty(K,\Psi)}_1(X)$, we let $f(X):=\sigma_X^{-1}(y)$.
        Function $f$ is defined $\mu_\kappa$-almost everywhere.

        \item It follows now that $j([f])=y$.
        This shows that $j$ is surjective.

        \item Thus, $j$ is the identity, which suffices for the conclusion of the claim.
    \end{prooff}

    \item For $X\in [\kappa]^\omega$, let us denote by
    \begin{parts}
        \item $\pi_X:\hh_X\to\hh_\Gamma$ the anti-collapse of $\Hull^{\hh_\Gamma}_1(X\cap\Theta_\Gamma)$,

        \item $\Xi_X$ the $\pi_X$-pullback strategy on $\hh_X$,\footnote{Note that $\sup\pi_X[o(\hh_X)]<\Theta_\Gamma$, so this pullback is a well-defined, complete strategy.}

        \item $K_X$ the union of all $J_1(N)$, where $(N,\Tau)$ is a mouse pair satisfying that $(\hh_X,\Xi_X)\unlhd^*(N,\Tau)$ and $\rho^-(N)=\rho_\omega(N)=o(\hh_X)$.\footnote{Observe that this is how we express ``$(N,\Tau)$ anticipates $(K,\Psi)$'' inside $\Mod$.}
    \end{parts}
    Observe that $(K_X : X \in [\kappa]^\omega)\in\Mod$.

    \item \claim For $\mu_\kappa$-almost every $X$, the set $K_X$ is the domain of the premouse $\kk_X$.
    \begin{prooff}
        \item Let $A$ consist of the sets $\ran(\pi_{(K',\Psi'),\infty})$, where $(K',\Psi')$ is an arbitrary non-dropping iterate of $(K,\Psi)$.
        Proposition \ref{535} implies that $A\in \mu_{K_\infty}$.

        \item It follows from Lemma \ref{562} that the set
        $$B:=\{Y \cap \kappa : Y \in A\}$$
        belongs to $\mu_\kappa$.

        \item\label{1917} The fact that every element of a premouse is $\Sigma_1$-definable from a finite sequence of ordinals yields that for all $X \in B$, the set $\Hull^{\M_\infty(K,\Psi)}_1(X)$ is the unique $Y \in A$ such that $X = Y \cap \kappa$.

        \item This shows that for all $X \in B$, it $(\kk_X, \Psi_X)$ is a non-dropping iterate of $(K, \Psi)$.

        \item Reasoning similarly to the line \ref{1917}, we conclude that for all $X \in B$, 
        $$\Hull^{\M_\infty(K,\Psi)}_1(X) \cap \hh_\Gamma = \Hull^{\hh_\Gamma}_1 (X \cap \Theta_\Gamma).$$

        \item The last two lines mean that for all $X \in B$, the mouse pair 
        $$(\kk_X, \Psi_X)|\delta^{\kk_X}$$
        is equal to $(\hh_X, \Xi_X)$ and that the domain of $\kk_X$ is equal to the union of $J_1(N)$, where $(N,\Tau)$ is an arbitrary mouse pair anticipating $(K, \Psi)$.
        For the second conclusion, we used Corollary \ref{nct-352}.

        \item What has been just said is simply another way of saying that for all $X \in B$, the domain of $\kk_X$ is equal to $K_X$.
        Since $B$ belongs to $\mu_\kappa$, the claim is proven.
    \end{prooff}

    \item Let
    $$j:\left(\prod_{\mu_\kappa}K_X\right)^\Mod\longrightarrow\prod_{\mu_\kappa}\lfloor\kk_X\rfloor$$
    be the canonical embedding.
    This embedding is $\Sigma_0$-elementary.\footnote{This comes down to the fact that $j$ commutes with rudimentary functions.}

    \item \claim $j$ is surjective.
    \begin{prooff}
        \item Let $y\in\prod_{\mu_\kappa}\lfloor\kk_X\rfloor$ be arbitrary and let us show that $y\in\ran(j)$.

        \item For $\mu_\kappa$-almost every $X$, the value $\sigma_X^{-1}(y)$ is a well-defined element of $K_X$.

        \item For such an $X$, we define $\hat{N}_{X,y}$ to be the $\in$-least bare premouse $\hat{N}$ coming from a mouse pair $(N,\Tau)$ anticipating $(K,\Psi)$ and satisfying that $\sigma_X^{-1}(y)\in J_1(\hat{N})$ (cf. Proposition \ref{nct-ordering by in}).

        \item We want to show that the family $(\hat{N}_{X,y})_X$ is almost everywhere equal to some family in $\Mod$.

        \item Let $\alpha:=o(\prod_{\mu_\kappa}\hat{N}_{X,y})< o(\M_\infty(K,\Psi))$.
        We have that $\alpha$ is represented in the ultraproduct $\prod_{\mu_\kappa}\lfloor\kk_X\rfloor$ by the function 
        $$h:[\kappa]^\omega\longrightarrow\Ord : X\mapsto \otp(\alpha\cap X).$$

        \item Observe that $h\in\Mod$ and that for $\mu_\kappa$-almost every $X$, $\hat{N}_{X,y}$ the $\in$-least bare premouse coming from a mouse pair $(N,\Tau)$ which anticipates $(K,\Psi)$ and satisfies that $o(\hat{N})=h(X)$.
        This shows that the function $X \mapsto \hat{N}_{X,y}$ is equal $\mu_\kappa$-almost everywhere to some function in $\Mod$.

        \item Since $\mu_\kappa$ is a supercompactness measure, there exist a Skolem term $t_y$ and a $p_y\in [\alpha]^{<\omega}$ such that for $\mu_\kappa$-almost every $X$,
        $$\sigma_X^{-1}(y)=t_y^{\hat{N}_{X,y}}(\sigma_X^{-1}(p_y)).$$

        \item For $\mu_\kappa$-almost every $X$,
        $$\sigma_X^{-1}(p_y)=\{\sigma_X^{-1}(\xi) : \xi\in p_y\}=\{\otp(\xi\cap X) : \xi\in p_y\}.$$
        Thus, the function $(\sigma_X^{-1}(p_y))_X$ belongs to $\Mod$ up to a negligible permutation.

        \item The previous three points show that the function $(\sigma_X^{-1}(y))_X$ belongs to $\Mod$ up to a negligible permutation.
        The corresponding equivalence class represents an element in $(\prod_{\mu_\kappa}K_X)^\Mod$ whose image via $j$ is equal to $y$.
        This concludes the verification that $j$ is surjective.
    \end{prooff}

    \item Since $j$ is $\Sigma_0$-elementary and surjective, it is the identity.
    Thus,
    $$K_\infty:=\lfloor\M_\infty(K,\Psi)\rfloor=\prod_{\mu_\kappa}\lfloor\kk_X\rfloor=\left(\prod_{\mu_\kappa}K_X\right)^\Mod\in\Mod.$$

    \item \claim $\Mod\models |K_\infty|\leq\kappa$
    \begin{proof}
        We work in $\Mod$.
        By the proof of \cite[Lemma 1.17]{schindler2010fine}, for every J-structure $M$, there exists a surjection $f:[o(M)]^{<\omega}\twoheadrightarrow M$ which is $\Sigma_1^M$, uniformly in $M$.
        Adding to this Proposition \ref{nct-ordering by in}, we get that for $\mu_\kappa$-almost every $X$, there exists a surjection
        $$f_X:[o(K_X)]^{<\omega}\twoheadrightarrow K_X,$$
        which is uniformly definable in $X$.
        By passing to the ultraproduct, we construct the surjection
        $$f:[o(K_\infty)]^{<\omega}\twoheadrightarrow K_\infty.$$
        This suffices for the conlusion.
    \end{proof}

    \item Since $\kappa^\omega\subseteq\Mod$, the last claim implies that $K_\infty^\omega\subseteq\Mod$.
    However, this means that $L(\Gamma,K_\infty^\omega)\subseteq\Mod$, which together with Proposition \ref{361} shows that $\ps(\R)^\Mod$ is strictly larger than $\Gamma$.
    This is a contradiction with \ref{nct-374}.
\end{prooff}

\newpage

\section{Synthesis}
\label{Synthesis}

We are now going to formulate precisely our construction of pointclasses and show that it reaches every model of $\AD^+ + \AD_\R$ in which there is no mouse pair with a measurable limit of Woodin cardinals.

\begin{declaration}
    We are assuming that $V\models\ZFC$.
    \qed
\end{declaration}

Let us remind the reader that the filter $\mu_X$ has been defined in Definition \ref{199} for all sets $X$.
Since we are assuming $\mathsf{AC}$, this filter is simply the club filter on $[X]^\omega$.
The construction to be defined can at most reach the minimal \textit{Chang type} pointclass.

\begin{definition}
    Let $\Delta$ be a closed pointclass.
    Then $\Delta$ is said to be of \intro{Chang type} iff for all $\kappa$, $\ps(\R)\cap L(\Delta^\omega,\kappa^\omega)[\mu_\Delta,\mu_\kappa]=\Delta$.
    \qed
\end{definition}

Our construction of pointclasses is essentially obtained by iterating the operator {$\Next$} to produce the next pointclass.
The definition of the operator that we propose cannot be realized if the current pointclass is of Chang type, so we will immediately exclude those pointclasses.

\begin{definition}
    Suppose that $\Delta$ is a closed pointclass which is \underline{not} of Chang type.
    Let $\kappa$ be the least ordinal such that
    $$\ps(\R)\cap L(\Delta^\omega,\kappa^\omega)[\mu_\Delta,\mu_\kappa]\supset\Delta.$$
    We define \intro{$\Next(\Delta)$} to be $\ps(\R)\cap L(\Delta^\omega,\kappa^\omega)[\mu_\Delta,\mu_\kappa]$.
    \qed
\end{definition}

When trying to reach by our construction the collection of all sets of reals of a certain determinacy model, it suffices to verify that a closed pointclass belonging to the model is not of Chang type and that the operator $\Next$ produces a pointclass that still belongs to the model in question.
The point is that $\Next$ can be applied to any non-Chang type pointclass and that it produces a strictly bigger pointclass by the definition.
This means that as long as the next pointclass is still inside the model, a simple recursive construction will relativize to the determinacy model and eventually reach all sets of reals.

\begin{proposition}\label{187}
    Suppose that $\Mod$ is an inner model containing all reals and satisfying $\AD^+ + \AD_\R+\mathsf{NMLW}$.
    Let $\Delta\in\Mod$ be a closed pointclass satisfying that $\Theta_\Delta<\Theta^\Mod$.
    Then $\Delta$ is \underline{not} of Chang type and $\Next(\Delta)\in\Mod$.
\end{proposition}
\begin{proof}
    Recall that for $\alpha<\Theta^\Mod$, we have that $\alpha^\omega\in\Mod$.
    This implies that also $\Delta^\omega\in\Mod$.
    Adding to this Proposition \ref{509}, we conclude that for all $\kappa<\Theta^\Mod$,
    $$L(\Delta^\omega,\kappa^\omega)[\mu_\Delta,\mu_\kappa]=L(\Delta^\omega,\kappa^\omega)[\mu_\Delta,\mu_\kappa]^\Mod.$$
    Hence, we will be done if we show that there exists $\kappa<\Theta^\Mod$ such that
    $$\ps(\R)\cap L(\Delta^\omega,\kappa^\omega)[\mu_\Delta,\mu_\kappa]\supset\Delta.$$
    
    If $\cof(\Theta_\Delta)=\omega$, then there exists a Wadge-cofinal sequence $(A_n : n<\omega)$ in $\Delta$.
    The joint of this sequence is a set of reals in $L(\Delta^\omega)$, which is not in $\Delta$.
    Hence, it remains to consider the case when $\cof(\Theta_\Delta)>\omega$.
    By Propositions \ref{509} and \ref{485}, $\omega_1$ is $\Delta$-supercompact in 
    $$L(\Delta)[\mu_\Delta]=L(\Delta)[\mu_\Delta]^\Mod.$$
    If $L(\Delta)\models\neg\AD_\R$, then Corollary \ref{565} implies that
    $$\Delta\subset\ps(\R)\cap L(\Delta)[\mu_\Delta].$$
    This means that the remaining case to verify is when $L(\Delta)\models\AD_\R$.
    In this case, Theorem \ref{1882} shows that there exists $\kappa<\Theta^\Mod$ such that
    $$\Delta\subset \ps(\R) \cap L(\Delta,\kappa^\omega)[\mu_\kappa].$$
\end{proof}

Let us now formulate the construction.
The first pointclass in the construction is the collection of sets of reals in $L(\R)$, provided that determinacy holds there.
If the determinacy fails in $L(\R)$, then the construction terminates immediately and is empty.
The construction also terminates if the current pointclass is of Change type.
Otherwise, we apply the operator $\Next$, which produces a bigger pointclass, but we take the bigger pointclass as a part of the construction only if it is a determinacy pointclass.
This means that the construction terminates as well if we ever reach a pointclass with a non-determined set.
At limit stages, we take the minimal closed pointclass extending all the previous ones, again terminating if a non-determined set is produced.

\begin{definition}
    We define \textcolor{Green}{$$(\De_\alpha : \alpha<\sh)$$}\index{$(\De_\alpha : \alpha<\sh)$} to be the longest sequence of closed pointclasses satisfying that for all $\alpha<\sh$,
    \begin{parts}
        \item $L(\De_\alpha)\models\AD^+$,


        \item if $\alpha=0$, then $\De_\alpha=\ps(\R)\cap L(\R)$,

        \item if $\alpha$ is a successor, then $\De_{\alpha-1}$ is \underline{not} of Chang type and $\De_\alpha=\Next(\De_{\alpha-1})$,

        \item if $\alpha$ is a limit, then, letting $\De_{<\alpha}:=\bigcup_{\xi<\alpha}\De_\xi$, it holds that
        $$\De_\alpha=\ps(\R)\cap L(\De_{<\alpha}).$$
        \qed
    \end{parts}
\end{definition}

The main result of the paper now reads as follows.

\begin{theorem}
    Suppose that $\Mod$ is an inner model containing all reals and satisfying $\AD^+ + \AD_\R+\mathsf{NMLW}$.
    Then there exists $\alpha<\sh$ such that 
    $$\ps(\R)\cap\Mod=\De_\alpha.$$
\end{theorem}
\begin{proof}
    This follows immediately from Proposition \ref{187}.
\end{proof}


\printindex
\newpage
\nocite{steel2023comparison}
\bibliographystyle{alpha}
\bibliography{lit.bib}
\end{document}